%% file: ProjectiveBundle4.tex
\documentclass[11pt,oneside]{amsart}
\usepackage[mathscr]{eucal}
\usepackage{amssymb, amsmath,array, amscd,verbatim}
\usepackage{graphics}

\newtheorem{thm}{Theorem}[section]

\newtheorem{prop}{Proposition}[section]

\newtheorem{remark}{Remark}[section]

\newtheorem{example}{Example}[section]

\newcommand{\F}{\ensuremath{\mathcal{F}}}
\newcommand{\C}{\ensuremath{\mathbb{C}}}
\newcommand{\cp}{\ensuremath{\mathbb{P}}}

\newcommand{\PSL}{\operatorname{PSL}(2,\C)}

\usepackage{color}

\usepackage[latin1]{inputenc}

\begin{document}

\title[Projectives structures and bundles]
{Projective structures and projective bundles over compact Riemann surfaces}

\author[Frank Loray and David Mar\'\i n]{Frank LORAY and David MAR\'IN}
\address{Frank LORAY (CNRS)
IRMAR, UFR de Math\'ematiques,
Universit\'e de Rennes 1, Campus de Beaulieu,
35042 Rennes Cedex (France)\newline
David MAR\'IN P\'EREZ,
Departament de Matem\'atiques,
        Universitat Aut\`onoma de Barcelona,
        Edifici Cc. Campus de Bellaterra,
        08193 Cerdanyola del Vall\`es (Spain)}
\email{frank.loray@univ-rennes1.fr\\ davidmp@mat.uab.es}
\urladdr{http://perso.univ-rennes1.fr/frank.loray 
http://mat.uab.es/~davidmp}

\date{\today}

\begin{abstract}A projective structure on a compact Riemann surface $C$
of genus $g$ is given by an atlas with transition functions in 
$\mathrm{PGL}(2,\C)$. Equivalently, a projective structure is given
by a $\mathbb P^1$-bundle over $C$ equipped with a section 
$\sigma$ and a foliation $\mathcal F$ which is both transversal
to the $\mathbb P^1$-fibers and the section $\sigma$.
From this latter geometric bundle picture, we survey on classical
problems and results on projective structures.
By the way, we will recall the analytic classification of $\mathbb P^1$-bundles.
We will give a complete description of projective (actually affine) structures
on the torus with an explicit versal family of foliated bundle picture.
\end{abstract}
\dedicatory{To Jos\'e Manuel AROCA  for his $60^{\text{th}}$ birthday}
\maketitle

\tableofcontents

\section{Projective structures}

\subsection{Definition and examples}

Denote by $\Sigma_g$ the orientable compact real surface of genus $g$.
A \textbf{projective structure} on $\Sigma_g$ is given by an atlas $\{(U_i,f_i)\}$ 
of coordinate charts (local homeomorphisms) $f_i:U_i\to\cp^1$ 
such that the transition functions 
$f_i=\varphi_{ij}\circ f_j$ are restrictions of Moebius transformations 
$\varphi_{ij}\in\mathrm{PGL}(2,\mathbb C)$.


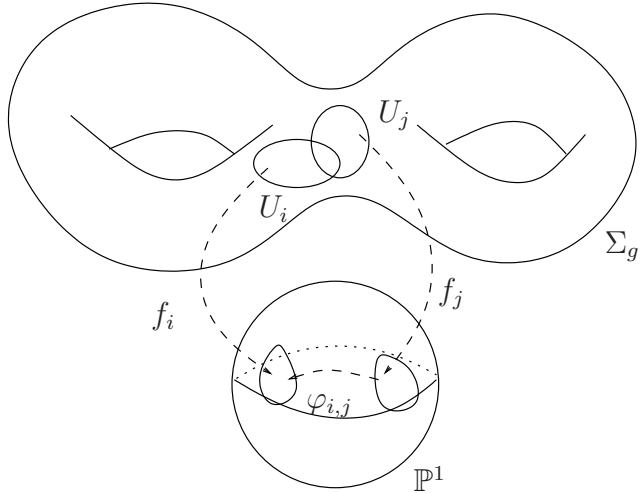
\begin{figure}[htbp]
\begin{center}

\input{Projective1.pstex_t}
 
\caption{Projective atlas}
\label{figure:1}
\end{center}
\end{figure}

There is a unique maximal atlas defining the projective structure above, 
obtained from the previous one by adding all charts $\{(U_i,\varphi\circ f_i)\}$ 
when $\varphi$ runs over $\mathrm{PGL}(2,\mathbb C)$.

A projective structure induces a \textbf{complex structure} on $\Sigma_g$,
just by pulling-back that of $\mathbb P^1$ by the projective charts.
We will denote by $C$ the corresponding Riemann surface (complex curve).

\begin{example}\label{E:UniversalCovering}(\textbf{The Universal cover}) \rm
Let $C$ be a compact Riemann surface having genus $g$
and consider its universal cover $\pi:U\to C$. By the Riemann Mapping Theorem,
we can assume that $U$ is either the Riemann sphere $\mathbb P^1$, 
or the complex plane $\mathbb C$ or the unit disk $\Delta$ depending
wether $g=0$, $1$ or $\ge2$. We inherit a representation of the fundamental group
$\rho:\pi_1(C)\to\text{Aut}(U)$ whose image $\Lambda$ 
is actually a subgroup of $\mathrm{PGL}(2,\mathbb C)$. All along the paper, by abuse of notation,
we will identify elements $\gamma\in\pi_1(C)$ with their image $\rho(\gamma)\in\mathrm{PGL}(2,\mathbb C)$.
The atlas defined on $C$ by all local determinations of
$\pi^{-1}:C\dashrightarrow \mathbb P^1$ 
defines a projective structure on $C$ compatible with the complex one.
Indeed, any two determinations of $\pi^{-1}$ differ by left composition 
with an element of $\Lambda$.
\end{example}

We thus see that any complex structure on $\Sigma_g$ is subjacent to a projective one.
In fact, for $g\ge1$, we will see that there are many projective structures compatible 
to a given complex one (see Theorem \ref{T:QuadraticDifferentials}). 
We will refer to the projective structure above as the \textbf{canonical projective structure}
of the Riemann surface $C$: it does not depend on the choice of the uniformization
of $U$. We now give other examples.


\begin{example}\label{E:Quotient}(\textbf{Quotients by Kleinian groups}) \rm
Let $\Lambda\subset\mathrm{PGL}(2,\mathbb C)$ be a subgroup acting properly, 
freely and discontinuously on some connected open subset $U\subset\mathbb P^1$.
Then, the quotient map $\pi:U\to C:=U/\Lambda$ induces a projective structure
on the quotient $C$, likely as in Example \ref{E:UniversalCovering}.
There are many such examples where $U$ is neither a disk, nor the plane.
For instance, \textbf{quasi-Fuchsian groups} are obtained as image of small perturbations 
of the representation $\rho$ of Example \ref{E:UniversalCovering};
following \cite{Weil60}, such perturbations keep acting discontinuously 
on some \textbf{quasi-disk}
(a topological disk whose boundary is a Jordan curve in $\mathbb P^1$). 
\end{example}

\begin{example}\label{E:SchottkyGroups}(\textbf{Schottky groups}) \rm
Pick $2g$ disjoint discs $\Delta_1^-,\ldots,\Delta_g^-$
and $\Delta_1^+,\ldots,\Delta_g^+$ in $\mathbb P^1$, $g\ge1$. For $i=1,\ldots,n$,
let $\varphi_i\in\mathrm{PGL}(2,\mathbb C)$ be a loxodromic map sending the disc $\Delta_i^-$
onto the complement $\mathbb P^1-\Delta_i^+$. 

\begin{figure}[htbp]
\begin{center}

\input{Projective2.pstex_t}
 
\caption{Schottky groups}
\label{figure:2}
\end{center}
\end{figure}
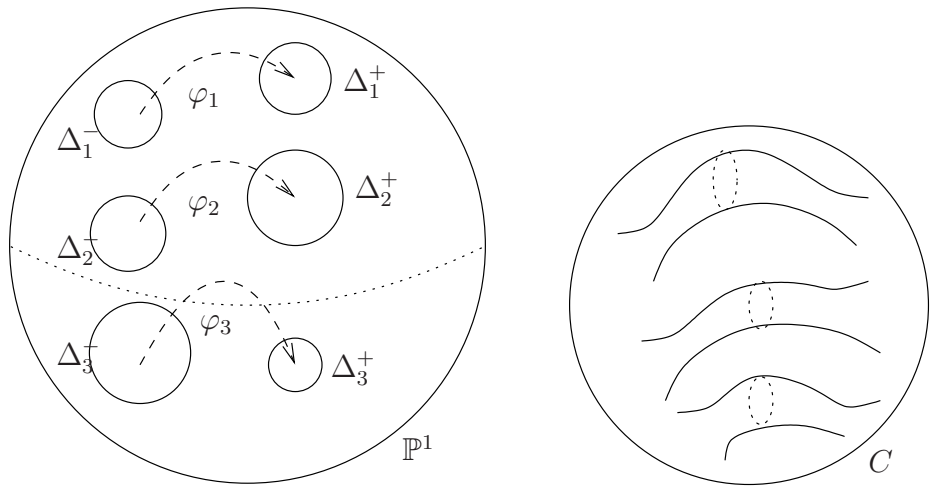

The group $\Lambda\subset\mathrm{PGL}(2,\mathbb C)$
generated by $\varphi_1,\ldots,\varphi_g$ acts properly, freely, and discontinuously
on the complement $U$ of some closed set contained inside the disks 
(a Cantor set whenever $g\ge2$).
The fundamental domain of this action on $U$ is given by the complement of the disks
and the quotient $C=U/\Lambda$ is obtained by gluing together the boundaries 
of $\Delta_i^+$ and $\Delta_i^-$ by means of $\varphi_i$, $i=1,\ldots,g$. 
Therefore, $C$ is a compact Riemann surface of genus $g$. This picture is clearly stable
under small deformation of the generators $\varphi_i$ and we thus obtain a 
complex $3g-3$ dimensional family of projective structures on the genus $g$ surface
$\Sigma_g$ (we have divided here by the action of $\mathrm{PGL}(2,\mathbb C)$ by conjugacy).
\end{example}


\subsection{Developping map and monodromy representation}\label{S:Monodromy}

Given a 
projective atlas and starting from any initial coordinate chart 
$(U_0,f_0)$, one can extend it analytically along any path $\gamma$ 
starting from $p_0\in U_0$.

Indeed, after covering $\gamma$ by finitely many projective coordinate charts, 
say $(U_0,f_0)$, $(U_1,f_1)$, ... ,$(U_n,f_n)$,
one can modify them
step by step in the following way.
First of all, since $f_0=\varphi_{01}\circ f_1$ on $U_0\cap U_1$, 
one can replace the chart $f_1$ by $\tilde f_1:=\varphi_{01}\circ f_1$
which is well-defined on $U_1$, extending $f_0$. Next, 
we replace $f_2$ by $\tilde f_2:=\varphi_{01}\circ\varphi_{12}\circ f_2$
which, on $U_1\cap U_2$, coincide with $\tilde f_1$.
Step by step, we finally arrive at the chart $\tilde f_n:=\varphi_{01}\circ\cdots\circ\varphi_{n-1\ n}\circ f_n$
which, by construction, is the analytic continuation of $f_0$ along $\gamma$.

Therefore, the local chart $(U_0,f_0)$
extends (after lifting on the universal covering) 
as a global submersion on the universal cover
$$f:U\to\mathbb P^1$$ 
which is called 
the \textbf{developping map} of the projective structure.
The developping map is moreover holomorphic with respect to the complex
structure subjacent to the projective one.
By construction, the monodromy of $f$ along loops takes the form
\begin{equation}\label{F:DeveloppingMap}
f(\gamma.u)=\varphi_\gamma\circ f,\ \ \ \varphi_\gamma\in\mathrm{PGL}(2,\mathbb C) 
\ \ \ \forall \gamma\in\pi_1(\Sigma_g,p_0)
\end{equation}
($u$ is the coordinate on $U$ and $\gamma.u$, the canonical action 
of $\pi_1(\Sigma_g,p_0)$ on $U$). In fact, 
$\varphi_\gamma$ is the composition of all transition maps $\varphi_{i,j}$
encoutered along $\gamma$ for a given finite covering of projective charts: 
with notations above, setting $(U_n,f_n)=(U_0,f_0)$, 
we have $\varphi_\gamma=\varphi_{01}\circ\cdots\circ\varphi_{n-1\ n}$.
It turns out that $\varphi_\gamma$ only depends on the homotopy class of $\gamma$
and we inherit a \textbf{monodromy representation}
\begin{equation}\label{F:Monodromy}
\rho:\pi_1(\Sigma_g,p_0)\to\mathrm{PGL}(2,\mathbb C)\ ;\ \gamma\mapsto\varphi_\gamma.
\end{equation}
The image $\Lambda$ of $\rho$ will be called monodromy group.
The developping map $f$ is defined by the projective structure up to the choice 
of the initial chart $(U_0,f_0)$ above: it is unique up to left composition
$\varphi\circ f$, $\varphi\in\mathrm{PGL}(2,\mathbb C)$. 
Therefore, the monodromy representation 
is defined by the projective structure up to conjugacy:
the monodromy of $\varphi\circ f$ is $\gamma\mapsto\varphi\circ\varphi_\gamma\circ\varphi^{-1}$.

Conversely, any global submersion $f:U\to\mathbb P^1$ on the universal covering
$\pi:U\to\Sigma_g$ satisfying (\ref{F:DeveloppingMap}) is the developping map
of a unique projective structure on $\Sigma_g$. We note that condition 
(\ref{F:DeveloppingMap}) forces the map $\gamma\to\varphi_\gamma$ to be a morphism.

\begin{example}\rm
The developping map of the canonical projective structure 
(see example \ref{E:UniversalCovering}) is the inclusion map 
$U\hookrightarrow\mathbb P^1$ of the universal cover of $C$.
More generally, when the projective structure is induced by a quotient map 
$\pi:U\to C=U/\Lambda$
like in example \ref{E:Quotient}, then the developping map $f$ is the universal
cover $\tilde U\to U$ of $U$ and the monodromy group is $\Lambda$.
In example \ref{E:SchottkyGroups}, the open set $U$ is not simply connected
(the complement of a Cantor set) and the developping map is a non trivial covering.
Thus the corresponding projective structure is not the canonical one. 
Similarly, the developping map of a quasi-Fuchsian group is the uniformization map
of the corresponding quasi-disk and is not trivial; the projective structure is 
neither the canonical one, nor of Schottky type.
\end{example}

\begin{example}\label{E:TheSphere}(\textbf{The Sphere}) \rm
Given a projective structure of the Riemann sphere $\mathbb P^1$,
we see that the developping map $f:\mathbb P^1\to\mathbb P^1$ is uniform
(no monodromy since $\pi_1(\mathbb P^1)$ is trivial). Therefore, $f$ is a
global holomorphic submersion (once we have fixed the complex structure)
and thus $f\in\mathrm{PGL}(2,\mathbb C)$. Consequently, the projective structure is the canonical one
and it is the unique projective structure on $\mathbb P^1$.
\end{example}

For similar reason, we remark that the monodromy group of a projective structure
on a surface of genus $g\ge1$ is never trivial.

\begin{example}\label{E:TheTorus}(\textbf{The Torus}) \rm
Let $\Lambda=\mathbb Z+\tau\mathbb Z$ be a lattice in $\mathbb C$, $\tau\in\mathbb H$,
and consider the elliptic curve $C:=\mathbb C/\Lambda$.
The monodromy of a projective structure on $C$ is abelian; 
therefore, after conjugacy, it is in one of the following abelian groups:
\begin{itemize}
\item the linear group $\{\varphi(z)=az\ ;\ a\in\mathbb C^*\}$,
\item the translation group $\{\varphi(z)=z+b\ ;\ b\in\mathbb C\}$,
\item the finite abelian dihedral group generated by $-z$ and $1/z$.
\end{itemize}
The canonical projective structure on $C$ has translation monodromy group $\Lambda$.
On the other hand, for any $c\in\mathbb C^*$ the map 
\begin{equation}\label{F:DeveloppingMapTorus}
f_c:\mathbb C\to\mathbb P^1\ ;\ u\mapsto\exp(c.u)
\end{equation}
is the developping map of a projective structure on $C$ whose monodromy is linear,
given by
\begin{equation}\label{F:LinearMonodromyTorus}
f_c(u+1)=e^{c}\cdot f(u)\ \ \ \text{and}\ \ \ f_c(u+\tau)=e^{c\tau}\cdot f(u).
\end{equation}
We inherit a $1$-parameter family of projective structures 
parametrized by $c\in\mathbb C^*$ (note that $f_0\equiv1$ is not a submersion).
We will see latter that this list is exhaustive. In particular, 
all projective structures on the torus are actually affine 
(transition maps in the affine group).
\end{example}

The projective structures listed in example \ref{E:TheTorus} 
are actually \textbf{affine structures}: the developping map takes values 
in $\mathbb C$ with affine monodromy.

\begin{thm}[Gunning \cite{Gunning67}]\label{T:GunningAffine} All projective structures on the elliptic curve
$\mathbb C/(\mathbb Z+\tau\mathbb Z)$, are actually affine and listed 
in example \ref{E:TheTorus} above.
There is no projective structure having affine monodromy 
on surfaces $\Sigma_g$ of genus $g\ge2$.
\end{thm}

In particular, the dihedral group is not the holonomy group 
of a projective structure on the torus.

\begin{proof}[Partial proof] Here, we only prove that the list 
of example \ref{E:TheTorus} exhausts all affine structures on 
compact Riemann surfaces. In example \ref{E:QuadraticDifferentialTorus}, 
we will see that there are no other
projective structure on tori than the affine ones.


Let $f:U\to\mathbb P^1$ be the developping map of a projective structure
with affine monodromy on the compact Riemann surface $C\not=\mathbb P^1$: 
$f$ is an holomorphic local homeomorphism satisfying 
$$f(\gamma\cdot u)=a_\gamma\cdot f(u)+b_\gamma,\ 
a_\gamma\in\mathbb C^*,\  b_\gamma\in\mathbb C,\
\ \ \forall\gamma\in\mathbb Z+\tau\mathbb Z.$$
Choose a holomorphic $1$-form $\omega_0$
on $C$ and write $\omega_0=\phi\cdot df$. Here, we identify $\omega_0$ 
with its lifting on the universal covering. Since $f$ is a local diffeomorphism,
$df$ has no zeroes and $\phi$ is a holomorphic on $U$,
vanishing exactly on zeroes of $\omega_0$ and poles of $f$. 
Moreover, the monodromy of $\phi$
is that of $df$, given by $\phi(\gamma\cdot u)=a_\gamma^{-1}\cdot\phi(u)$.
Therefore, the meromorphic $1$-form $\omega_1={d\phi\over\phi}$
has no monodromy: it defines a meromorphic $1$-form on $C$
having only simple poles, located at the zeroes of $\omega_0$ and poles of $f$,
with residue $+1$. Following Residue Theorem, $\omega_1$ has actually no poles:
$f$ is holomorphic, $\omega_0$ does not vanish and thus genus $g=1$. 
This proves the second assertion of the statement.

Now, assume $g=1$, $U=\mathbb C$ and $C=\mathbb C/(\mathbb Z+\tau\mathbb Z)$.
The $1$-form $\omega_1$ above is holomorphic and thus takes the form $\omega_1=-c\cdot du$
for some constant $c\in\mathbb C$. In other words, we have
$f''/f'=c$ and we obtain after integration
\begin{itemize}
\item $f(u)=a\cdot e^{cu}+b$ when $c\not=0$,
\item $f(u)=a\cdot u+b$ when $c=0$
\end{itemize}
for constants $a\in\mathbb C^*$ and $b\in\mathbb C$. After left composition
by an affine map, which does not affect the affine structure, we can set
$a=1$ and $b=0$ and $f$ 
belongs to the list of example \ref{E:TheTorus}.
\end{proof}

\begin{remark}\rm We see from the proof that the projective structures on 
$\mathbb C/(\mathbb Z+\tau\mathbb Z)$ are naturally parametrized by $\mathbb C$,
namely the constant map $\phi=f''/f'\equiv c$, which is not clear from 
the description of example \ref{E:TheTorus} (we see $\mathbb C^*$ plus one point
). 
One can recover this by choosing conveniently the integration constants $a$ and 
$b$ 
in the proof above. Indeed, consider the alternate family of developping maps
given by 
\begin{equation}\label{F:UniversalDeveloppingMapTorus}
F:\mathbb C^2\to\mathbb C\ ;\ (c,u)\mapsto 
f_c(u):=\left\{\begin{array}{cc}{e^{cu}-1\over c},& c\not=0\\ u,&c=0\end{array}\right.
\end{equation}
The map $F$ is clearly holomorphic on $\mathbb C^2$ and 
makes the developping maps $f_c$ into an holomorphic family parametrized 
by $c\in\mathbb C$. Moreover, the corresponding holonomy representations
are given by
$$f_c(u+\gamma)=\left\{\begin{array}{cc}e^{c\gamma}f_c(u)+{e^{c\gamma}-1\over c}
,& 
c\not=0\\ u+\gamma,&c=0\end{array}\right.\ \ \ \forall\gamma\in\mathbb Z+\tau\mathbb Z$$
and we see the affine motions with common fixed point $-1/c$
converging to translations while $c\to 0$.
\end{remark}


\begin{remark}\rm When we set $g=1$ in example \ref{E:SchottkyGroups},
we have $U=\mathbb C^*$ and $\Lambda$ is generated by a single map 
$\varphi(z)=e^{2i\pi\lambda}z$. The quotient $C=U/\Lambda$ is the
elliptic curve with lattice $\mathbb Z+\lambda\mathbb Z$. The complex
structure varies with $\lambda$ and very few projective structures
on the torus are obtained by this way. 
In fact, we see in example 
\ref{E:TheTorus} that, for generic values of $c$, the monodromy group 
of the corresponding projective structure is not discrete
($c$ is not $\mathbb Z$-commensurable with $1$ and $\tau$).
\end{remark}

\subsection{Quadratic differentials}

In order to generalize the arguments involved in the proof 
of Theorem \ref{T:GunningAffine}
for genus $g\ge2$ Riemann surfaces, we have to replace $f''/f'$
by the \textbf{Schwartzian derivative} of $f$
\begin{equation}\label{F:SchwartzianDerivative}
\mathcal S(f):=\left({f''\over f'}\right)'-{1\over 2}\left({f''\over f'}\right)^2.
\end{equation}
Recall that, for any holomorphic functions $f$ and $g$, we have
\begin{equation}\label{F:SchwartzianComposition}
\mathcal S(f\circ g)=\mathcal S(f)\circ g\cdot (g')^2+\mathcal S(g).
\end{equation}
Given a projective structure on a Riemann surface $C$, consider
the Schwartzian derivative of the corresponding developping map
$\phi:=\mathcal S(f)$. For any $\gamma\in\Lambda=\pi_1(C)$, 
we deduce from property (\ref{F:DeveloppingMap}) of $f$ that
$$\phi\circ\gamma\cdot(\gamma')^2=\mathcal S(f\circ\gamma)
=\mathcal S(\varphi_\gamma\circ f)=\phi.$$
In other words, the \textbf{quadratic differential} $\omega=\phi(u)\cdot du^2$
is invariant under $\Lambda$ and gives rise to a quadratic differential
on the Riemann surface $C$. We note that $\omega$ is holomorphic.
Indeed, outside the poles of $f$, $\phi_1:=f'$ is not vanishing, 
thus $\phi_2:=\phi_1'/\phi_1$ is holomorphic and $\phi=\phi_2'-(\phi_2)^2/2$
well. On the other hand, at a pole of $f$, one can replace $f$
for instance by $1/f$, which is not relevant for the Schwartzian derivative, 
and go back to the previous argument. By this way, we canonically associate 
to any projective structure on $C$ an holomorphic quadratic differential
$\omega$ on $C$, i.e. a global section of $K_C^{\otimes 2}$, 
where $K_C$ is the canonical line bundle over $C$.

Conversely, given any holomorphic quadratic differential $\omega=\phi(u)\cdot du^2$
on the Riemann surface $C$, one can solve locally the differential equation
$\mathcal S(f)=\phi$ in $f$ and recover the coordinate charts of a projective
structure on $C$ (compatible with the complex one): the fact is that any two
(local) solutions of $\mathcal S(f)=\phi$ differ by left composition by 
a Moebius transformation.

\begin{example}\label{E:QuadraticDifferentialTorus}\rm
In genus $1$ case, any holomorphic quadratic differential
takes the form $\omega=c\cdot du^2$ for a constant $c\in\mathbb C$
($K_C^{\otimes 2}=K_C$ is still the trivial bundle). 
In fact, $\omega=\tilde{\omega}^{2}$, where $\tilde{\omega}=\sqrt{c}du$.
On the other hand, any solution of $f''/f'=\tilde c$ gives rise to a solution 
of $\mathcal S(f)=-\tilde c^2/2=c$; therefore, the projective
structure defined by $\omega$ is actually subjacent to the
affine structure defined by $\tilde c$. This concludes the proof
of Theorem \ref{T:GunningAffine}. We note that the space of affine structures 
forms a two fold covering of the space of projective structures
(the choice of the square root $\tilde c$). Of course, this comes from the fact that
the $2$ affine structures given by $f_c$ and $1/f_c$ (with notations of 
example \ref{E:TheTorus}) do not define distinct projective structures.
\end{example}

For genus $g\ge2$ Riemann surfaces, the dimension of $H^0(C,K_C^{\otimes 2})$
can be computed by Riemann-Roch Formula, and we obtain

\begin{thm}[Gunning \cite{Gunning67}]\label{T:QuadraticDifferentials} 
The set of projective structures 
on a complex Riemann surface $C$ of genus $g\ge2$ is parametrized
by the $3g-3$-dimensional complex vector space $H^0(C,K_C^{\otimes 2})$.
\end{thm}

In this vector space, $0$ stands for the canonical structure of example
\ref{E:UniversalCovering}.


\subsection{The monodromy mapping}

A natural question arising while studying projective structures is to understand,
for a given surface $\Sigma_g$, the nature of the Monodromy (or Riemann-Hilbert) Mapping
$$\mathcal{P}_{g}\longrightarrow\mathcal{R}_{g}.$$
On the left-hand side, $\mathcal P_g$ denotes the set of all projective structure on $\Sigma_g$
up to isomorphism; on the right-hand side, 
$\mathcal R_g$
is the set of representations
of the fundamental group in $\mathrm{PGL}(2,\mathbb C)$ up to conjugacy:
$$\mathcal R_g=
\text{Hom}(\pi_1(\Sigma_g),\mathrm{PGL}(2,\mathbb C))/_{\mathrm{PGL}(2,\mathbb C)}.$$ 

Let us first consider the genus $g=1$ case. 
From Gunning's Theorem \ref{T:GunningAffine},
the left-hand side can be viewed as a $\mathbb C$-bundle over the modular orbifold
$\mathbb H/\text{PSL}(2,\mathbb Z)$ where $\mathbb H$ 
denotes the upper-half plane whose fiber at a given complex structure
is the affine line of holomorphic differentials. 
Nevertheless, to avoid dealing with orbifold points, we prefer
to deal with the parametrization of affine structures 
by $\mathbb H\times\mathbb C$
given by the map
$$(\tau,c)\mapsto (C,\omega)\ \text{where}\ \left\{
\begin{array}{cc} C=&\mathbb C/(\mathbb Z+\tau\mathbb Z)\\
\omega=&c\cdot du\end{array}\right.$$
Here, the base $\mathbb H$ is the space of marked complex structures
on the torus, up to isomorphism, and the fiber over $\tau$ is the 
affine line of differentials $\mathbb C\cdot du$, 
$u$ the variable of $\mathbb H$.
Since all projective structures are actually
affine, we can replace$\mathcal{R}_{1}$ by
$\mathcal A_1:=\text{Hom}(\pi_1(C),\text{Aff}(\mathbb C))/_{\text{Aff}(\mathbb C)}$
where 
$$\text{Aff}(\mathbb C):=\{\varphi(z)=az+b,\ a\in\mathbb C^*,\ b\in\mathbb C\}$$
is the group of affine transformations. Once we have fixed generators $1$ and $\tau$
for the fundamental group of $C=\mathbb C/(\mathbb Z+\tau\mathbb Z)$,
the set $\text{Hom}(\pi_1(C),\text{Aff}(\mathbb C))$ identifies with the complex
$3$-dimensional subvariety 
$$\{(a_1z+b_1,a_\tau z+b_\tau)\ ;\ (a_1-1)b_\tau=(a_\tau-1)b_1\}\subset(\mathbb C^*\times\mathbb C)^2$$
(here, we see the condition for the commutativity). 
Since any linear representation
does not occur as monodromy representation
of an affine structure on the torus, we consider the quotient $\mathcal B_{1}\subset\mathcal A_1$ 
of the complement of $b_{1}=b_{\tau}=0$ in this variety.
It is easy to see that $\mathcal B_1$ can be
identified with the $2$-dimensional complex manifold
$$\mathcal B_1=\{(a_1,a_\tau,[b_1:b_\tau])\ ;\ (a_1-1)b_\tau=(a_\tau-1)b_1\}\subset
\mathbb C^*\times \mathbb C^*\times \mathbb P^1$$
where $[z:w]$ denotes homogeneous coordinates on $\mathbb P^1$.

The projection $\mathcal B_1\to\mathbb C^*\times \mathbb C^*$ is just the blow-up
of the point $(1,1)$ and the exceptional divisor consists in euclidean representation.
Finally, the monodromy map is described by
$$\mathbb H\times\mathbb C\to\mathcal B_1\ ;\ (\tau,c)\mapsto
\left\{\begin{array}{cc}(e^c,e^{c\tau},[{e^{c}-1\over c}:{e^{c\tau}-1\over c}]),
& c\not=0\\ 
(1,1,[1:\tau]),&c=0\end{array}\right.$$


Looking at the differential of the Monodromy Map above, 
we see that it has always rank $2$ and {\it the Monodromy Map is 
an holomorphic local diffeomorphism}; it is moreover {\it injective 
and proper in restriction to each fiber $\tau\times\mathbb C$}.
Its image is the complement of the real torus 
$\mathbb S^1\times\mathbb S^1\subset\mathbb C^*\times \mathbb C^*$,
or we better should say, its lifting on $\mathcal B_1$:
the complement of $\mathbb P^1(\mathbb R)$ in restriction to the exceptional divisor.

{\bf But} the Monodromy Map {\it is neither injective,
nor  a covering map onto its image}: for instance, 
for any $\tau,\tau'\in\mathbb H$, $\tau'\not=\tau$, 
and for any $(m,n)\in\mathbb Z^2-\{(0,0)\}$, 
the two affine structures
$$(\tau,2i\pi{m\tau'+n\over\tau'-\tau})\ \ \ \text{and}\ \ \ 
(\tau',2i\pi{m\tau+n\over\tau'-\tau})$$
have the same monodromy representation. In particular, 
the injectivity is violated for arbitrarily close complex structures.
On the other hand, the monodromy of the canonical
structure $(\tau,0)$ occurs only for this structure.

Consider now the genus $g\ge2$ case. 
It follows from Gunning's Theorem \ref{T:QuadraticDifferentials} above 
that the set $\mathcal P_g$ 
of projective structures on the genus $g\ge2$ surface $\Sigma_g$ 
can be viewed as a complex $6g-6$-dimensional space. Indeed, if we denote 
by $\mathcal T_g$ the Teichm\"uller space of complex marked structures
on $\Sigma_g$ viewed as an open subset of $\mathbb C^{3g-3}$, 
then $\mathcal P_g$ is parametrized by the rank $3g-3$-vector bundle
$\tilde{\mathcal P}_g$ over $\mathcal T_g$ whose fiber over 
a given complex structure $C$ is the space of quadratic differentials
$H^0(C,K_C^{\otimes 2})$.

By Theorem \ref{T:GunningAffine}, the monodromy representation cannot be affine in the case $g\ge 2$. The image of $\mathcal{P}_{g}$ by the Monodromy Map  is thus included in the subset of 
\textbf{irreducible} representations
$$\mathcal R_{g}^{\mathrm{irr}}:=\mathcal R_{g}-\mathcal A_{g}$$
where $\mathcal A_{g}=\mathrm{Hom}(\pi_{1}(\Sigma_{g}),\mathrm{Aff}(\C))/_{\mathrm{PGL}(2,\mathbb C)}$ is the set of affine representations {\bf up to $\mathrm{PGL}(2,\mathbb C)$-conjugacy}.
One can check (see \cite{Gunning71}) that $\mathcal R_{g}^{\mathrm{irr}}$ forms a non-singular complex manifold 
of dimension $6g-6$.
Thus, the Monodromy Map can locally be described as a holomorphic map between open subsets of  $\C^{6g-6}$ and the following result makes sense
(see proof in section \ref{S:BundlePicture}).

\begin{thm}[Hejhal \cite{Hejhal75,Earle78,Hubbard78}]\label{T:LocalDiffeomorphism} 
The Monodromy Map is a local diffeomorphism.
\end{thm}

In \cite{Kawai96}, it is moreover proved that the Monodromy Map is symplectic
with respect to symplectic structures that can be respectively canonically defined
on both spaces (see \cite{Goldman84}).

The restriction of the Monodromy Map to each fiber $H^{0}(C,K_{C}^{\otimes 2})$ of $\tilde{\mathcal{P}}_{g}$ over $C\in\mathcal{T}_{g}$ is injective. In other words, we have the following result whose proof will be given in section \ref{sec2.2}.

\begin{thm}[Poincar\'e \cite{Poincare84}]\label{T:Poincare} 
Given a compact Riemann surface $C$, 
any two projective structures are the same if, and only if,
they have the same monodromy representation (up to $\mathrm{PGL}(2,\mathbb C)$).
\end{thm}

It is clear that the Monodromy Map is not surjective. First of all,
by Theorem \ref{T:GunningAffine},
its image is contained in $\mathcal{R}_{g}^{\mathrm{irr}}\subset\mathcal{R}_{g}$. 
On the other hand, the space of representations 
$\text{Hom}(\pi_1(C),\mathrm{PGL}(2,\mathbb C))$ falls into $2$ connected components,
namely the component of those that can be lifted as 
$\text{Hom}(\pi_1(C),\text{SL}(2,\mathbb C))$ and the other ones.
Since 
the Monodromy Map is continuous
(actually holomorphic) and since the monodromy of canonical
projective structures can be lifted to $\text{SL}(2,\mathbb R)$,  it becomes clear
that the image of the Monodromy Map will be in the former component.
Finally, notice that the monodromy
representation cannot be in $\text{PSU}(2,\mathbb C)$, i.e. conjugated to a group of rotations of the sphere,
otherwise we could pull-back the invariant spherical metric of $\mathbb P^1$
by the developping map, giving rise to a curvature $+1$ metric on the 
surface, impossible except in the trivial case $g=0$ (see Example \ref{E:TheSphere}).
The main result in the field, which has been conjectural for decenies, is the following.

\begin{thm}[Gallo-Kapovich-Marden \cite{GKM00}]\label{T:GalloKapovichMarden} 
Consider the genus $g$ surface $\Sigma_g$, $g\ge2$. 
An homomorphism $\rho\in\mathrm{Hom}(\pi_1(\Sigma_g),\mathrm{PGL}(2,\mathbb C))$
is the monodromy representation of a projective structure on the $\Sigma_g$
if, and only if, $\rho$ can be lifted as 
$\tilde\rho\in\mathrm{Hom}(\pi_1(C),\mathrm{SL}(2,\mathbb C))$ and the image of $\rho$
is, up to $\mathrm{PGL}(2,\mathbb C)$-conjugacy, neither in the affine group $\mathrm{Aff}(\mathbb C)$,
nor in the rotation group $\mathrm{PSU}(2,\mathbb C)$.
\end{thm}

\subsection{The fibre bundle picture}\label{S:BundlePicture}

Let $f:U\to\mathbb P^1$ be the developping map of a projective structure 
on $C$ (here we fix the underlying complex structure) and consider
its graph $\{(u,f(u))\ ;\ u\in U\}\subset U\times\mathbb P^1$.
The fundamental group $\pi_1(C)$ acts on the product $U\times\mathbb P^1$
in as follows: for any $\gamma\in\pi_1(C)$, set
$$\gamma: (u,y)\mapsto (\gamma\cdot u,\varphi_\gamma(y))$$
where $u\mapsto\gamma\cdot u$ is the canonical action of $\pi_1(C)$
on the universal cover and $\varphi_\gamma$ is the 
monodromy of the projective structure along $\gamma$.
This action of $\pi_1(C)$ is proper, free and discontinuous 
since its projection on $U$ is so. By consequence, we can consider 
the quotient:
$$P:=U\times\mathbb P^1/_{\pi_1(C)}.$$
The projection $U\times\mathbb P^1\to C$ defined by 
$(u,y)\mapsto \pi(u)$, where $\pi:U\to C$ is the universal cover,
is preserved by the action and induces a global submersion
$$\pi:P\to C$$ 
making {\it $P$ into a $\mathbb P^1$-bundle over $C$.
The graph of $f$ also is invariant under the action
(consequence of (\ref{F:DeveloppingMap})) thus defining a section }
$$\sigma:C\to P.$$
Finally, the horizontal foliation defined by $\{y=\text{constant}\}$
is also preserved and defines a foliation $\mathcal F$ transversal 
to all $\mathbb P^1$-fibres on $P$.
Since the developping map $f$ is regular, its graph is transversal to the horizontal
foliation and $\sigma$ is transversal to $\mathcal F$. 
In this situation, we say that the $\mathbb P^{1}$-bundle $P$ is flat.
The triple $(\pi:P\to C,\mathcal F,\sigma)$ is well-defined
by the projective structure up to analytic isomorphism of $\mathbb P^1$-bundles.

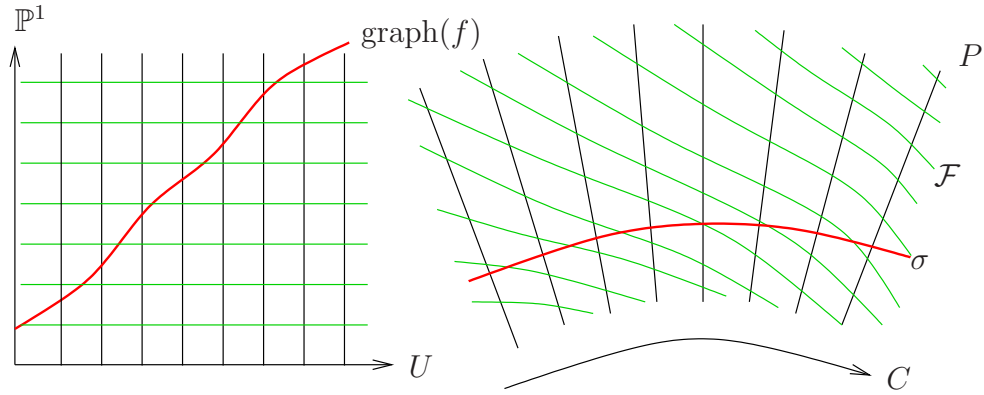
\begin{figure}[htbp]
\begin{center}

\input{Projective3bis.pstex_t}
 
\caption{From projective structure to bundle picture}
\label{figure:3}
\end{center}
\end{figure}

Conversely, given a $\mathbb P^1$-bundle $\pi:P\to C$, a foliation $\mathcal F$ 
on $P$ transversal
to $\pi$ and a section $\sigma:C\to P$ transversal to $\mathcal F$, then the (unique) projective structure
on $\mathbb P^1$-fibres can be transported, transversely to the foliation $\mathcal F$, inducing a projective
structure on the section $\sigma(C)$, and thus on its $\pi$-projection $C$.

\begin{figure}[htbp]
\begin{center}

\input{Projective4.pstex_t}
 
\caption{From bundle picture to projective structure}
\label{figure:4}
\end{center}
\end{figure}
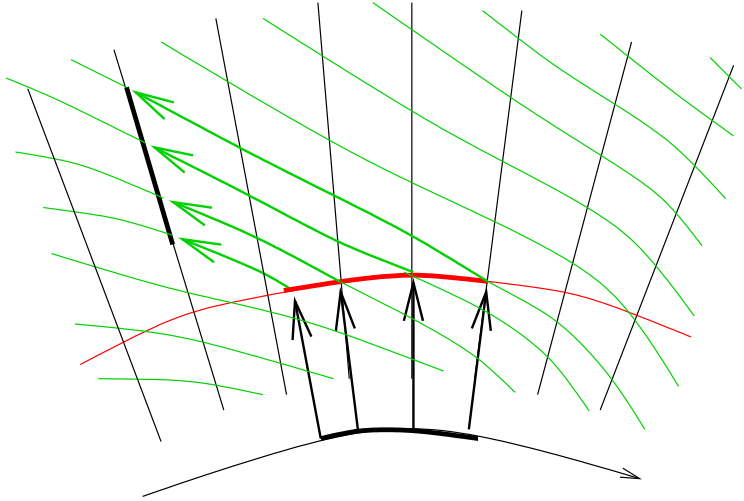

In the recent terminologoly of \cite{BeilinsonDrinfeld},
such triple $(\pi:P\to C,\mathcal F,\sigma)$ are called 
$\mathfrak{sl}(2,\mathbb C)$-opers.

\begin{remark}\label{rm13}\rm Given an homomorphism $\rho\in\text{Hom}(\pi_1(C),\text{PGL}(2,\mathbb C))$,
one can at least construct the pair $(\pi:P\to C,\mathcal F)$ as above. This foliated surface is called the
{\bf suspension} of the representation $\rho$, also known as the {\bf flat
$\mathbb P^1$-bundle} associated to $\rho$. Conversely, consider a flat 
$\mathbb P^1$-bundle, i.e.
a pair $(\pi:P\to C,\mathcal F)$ where $\pi:P\to C$ is a $\mathbb P^1$-bundle
and $\mathcal F$ is a foliation transversal to $\pi$. Then one can associate to it
a representation $\rho$ in the following way.

Over any sufficiently small open subset $U_i\subset C$, one can construct
a trivializing coordinate $F_i:\pi^{-1}(U_i)\to \mathbb P^1$ for the flat bundle,
that is to say inducing an isomorphism in restriction to each fibre and
such that the level curves $F_i^{-1}(y_0)$ are local leaves of the
foliation $\mathcal F$. In fact, $F_i$ is uniquely determined after choosing 
the local $\mathcal F$-invariant sections $\sigma_0,\sigma_1,\sigma_\infty:U_i\to P$
along which $F_i$ takes values $0$, $1$ and $\infty$ respectively. 
Such flat coordinate is well defined up to left-composition by a Moebius transformation;
likely as in section \ref{S:Monodromy}, after fixing a flat local coordinate $F$
over some neighborhood of the base point $x_0\in C$,
we inherit a monodromy representation
$\rho:\pi_1(C,x_0)\to\mathrm{PGL}(2,\mathbb C)$ 
where the analytic continuation of $F$ along any loop $\gamma$
satisfies $F(\gamma\cdot u)=\rho(\gamma)\circ F(u)$.

It turns out that any flat $\mathbb P^1$-bundle is isomorphic to the suspension
of its monodromy representation just defined.  
In fact, {\it any two flat $\mathbb P^1$-bundles are isomorphic if, and only if,
they have the same monodromy representation up to $\mathrm{PGL}(2,\mathbb C)$
conjugacy}. Indeed, let $(\pi:P\to C,\mathcal F)$ and $(\pi':P'\to C,\mathcal F')$ 
be flat $\mathbb P^1$-bundles having flat coordinates $F$ and 
$F'$ over $U_0\subset C$ giving rise to the same monodromy representation;
then the local isomorphism $\Phi:\pi^{-1}(U_0)\to\pi'^{-1}(U_0)$ sending 
any point $p$ to the unique point $p'$ satisfying $(\pi(p),F(p))=(\pi'(p'),F'(p'))$
extends uniformly as a global isomorphism of flat $\mathbb P^1$-bundles
$\Phi:P\to P'$, i.e. conjugating $\mathcal F$ to $\mathcal F'$ and satisfying
$\pi'\circ\Phi=\pi$.
 \end{remark}


\begin{proof}[Proof of Hejal's Theorem \ref{T:LocalDiffeomorphism}]
In fact, since the Monodromy Map is clearly holomorphic,
it is enough to prove that it is locally bijective.

Let $(\pi:P\to C,\mathcal F,\sigma)$ be the triple associated to a
projective structure having monodromy representation 
$\rho\in\text{Hom}(\pi_1(\Sigma_g),\text{PGL}(2,\mathbb C))$.
For any perturbation $\rho'\in\text{Hom}(\pi_1(\Sigma_g),\text{PGL}(2,\mathbb C))$
of $\rho$, the corresponding suspension $(\pi':P'\to C,\mathcal F')$ is 
close to the foliated bundle $(\pi:P\to C,\mathcal F)$; if the perturbation
is small enough, one can find a real $C^\infty$ section $\sigma':C\to P'$
close to $\sigma:C\to P$ and still transversal to $\mathcal F'$ 
(all of this makes sense and can be checked on the neighborhood
of a fundamental domain of the universal cover $U\times\mathbb P^1$).
The foliation $\mathcal F'$ still induces a projective structure
on the real surface $\sigma'(C)$ that, by construction, has the required
monodromy. This proves the surjectivity.

Let $(\pi:P\to C,\mathcal F,\sigma)$ be the triple associated to a
projective structure $\mathcal P$ and consider another projective structure 
$\mathcal P'$ close to this $\mathcal P$ having the same monodromy
representation. The fibre bundle construction
can be done in the real $C^\infty$ setting so that one can associate
to
$\mathcal P'$ a triple $(\pi:P\to C,\mathcal F,\sigma')$ where $C$
is still the complex curve attached to $\mathcal P$ and $\sigma':C\to P$
is now a real $C^\infty$ section transversal to $\mathcal F$;
we note that the pair $(\pi:P\to C,\mathcal F)$ is the same for $\mathcal P$
and $\mathcal P'$ since they have the same monodromy representation.
If $\mathcal P'$ is close enough to $\mathcal P$, say in the $C^\infty$ category, 
then $\sigma'$ is close to $\sigma$; one can therefore 
unambiguously define a $C^\infty$ diffeomorphism $\phi:\sigma'(C)\to\sigma(C)$
by following the leaves of the foliation from one section to the other one.
By construction, the projective structures induced by $\mathcal F$
on both sections are conjugated by $\phi$.
The diffeomorphism $\pi_*\phi:=\pi\circ\phi\circ\sigma'$ actually
integrates the quasi-conformal structure induced by $\mathcal P'$
on $C$; it is close to the identity.
\end{proof}

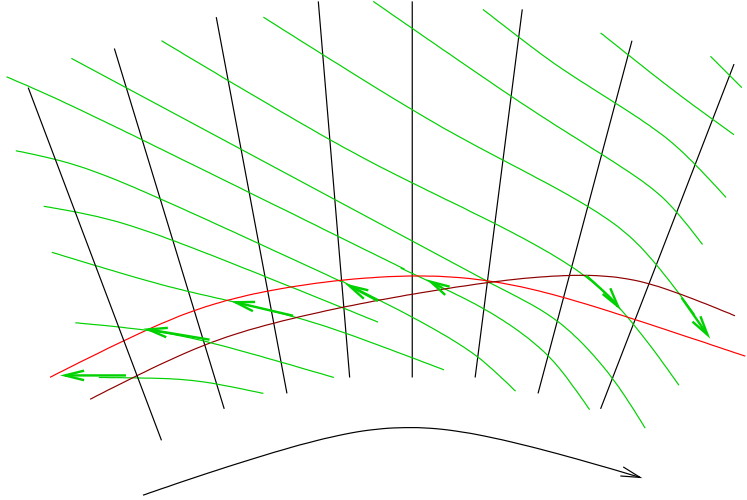
\begin{figure}[htbp]
\begin{center}

\input{Projective4bis.pstex_t}
 
\caption{Local injectivity of the Monodromy Map}
\label{figure:4bis}
\end{center}
\end{figure}

\begin{remark}\rm Since the Monodromy Map is not globally injective,
the injectivity argument of the previous proof cannot be carried out
for sections $\sigma$ and $\sigma'$ that are not closed enough:
the set of $C^\infty$ sections transversal to $\mathcal F$ 
may have infinitely many connected components as it so happens
in the case of affine structures on the torus.
Similarly, the surjectivity argument of the proof cannot be globalized:
when the monodromy representation $\rho'$ eventually becomes
reducible for instance, there does not exist $C^\infty$ section 
transversal to $\mathcal F$ anymore.
Following Theorem \ref{T:GalloKapovichMarden}, the existence 
of a $C^\infty$ section transversal to $\mathcal F$ is possible
if, and only if, $\mathcal F$ is the suspension of a non elementary 
representation $\rho$ (lifting to $\text{SL}(2,\mathbb C)$) ! 
From this point of view, Theorem \ref{T:GalloKapovichMarden}
looks like a very subtle transversality result.
\end{remark}


\section{$\mathbb P^1$-bundles and Riccati foliations}

Motivated by the fibre bundle picture of section \ref{S:BundlePicture},
we developp here the study of Riccati foliations on $\mathbb P^1$-bundles
over compact Riemann surfaces.

\subsection{Classification of $\mathbb P^1$-bundles}

Let $\pi:P\to C$ be a $\mathbb P^1$-bundle over a compact Riemann surface $C$:
$P$ is a smooth surface and the fibers of $\pi$ are rational, isomorphic to $\mathbb P^1$.
We also say that $P$ is a {\bf ruled surface}. Another $\mathbb P^1$-bundle $\pi':P'\to C$ 
is {\bf analytically equivalent} to the previous one if there is an holomorphic 
diffeomorphism
$\phi:P\to P'$ such that $\pi'\circ\phi=\pi$.
We recall some basic facts (see \cite{Hartshorne,Maruyama70}). 

On open charts $u_i:U_i\to\mathbb C$ on $C$, the bundle becomes analytically trivial (see \cite{Forster}):
we have holomorphic diffeomorphisms (trivializing coordinates)
$$\phi_i:\pi^{-1}(U_i)\to U_i\times\mathbb P^1\ ;\ p\mapsto (\pi(p),\varphi_i(p)).$$ 
On overlapping charts $U_i\cap U_j$, the transition maps
take the form $\phi_i=\phi_{i,j}\circ\phi_j$ where
$\phi_{i,j}(u,y)=(u,\varphi_{i,j}(u,y))$ and 
$$\varphi_{i,j}\in\text{PGL}(2,\mathcal O(U_{i,j})).$$ 
The $\mathbb P^1$-bundle is equivalently defined by the collection 
$$(\varphi_{i,j})_{i,j}\in H^1(C,\text{PGL}(2,\mathcal O)).$$
By lifting conveniently the transition maps into $H^1(C,\text{GL}(2,\mathcal O))
$,
we may view a $\mathbb P^1$-bundle as the projectivization $P=\mathbb PV$
of a {\bf rank $2$ vector bundle} $V$ over $C$. Moreover, another vector bundle 
$V'$
will give rise to the same $\mathbb P^1$-bundle if, and only if, 
$V'=L\otimes V$ for a line bundle $L$ over $C$. 
The classification of $\mathbb P^1$-bundles is thus equivalent
to the classification of rank $2$ vector bundles up to tensor product by a line 
bundle.

From the {\bf topological} point of view, due to the fact that $\pi_1(\text{PGL}(2,\mathbb
C))=\mathbb Z/2\mathbb Z$, there are exactly $2$ distinct $\mathbb S^2$-bundles
over a compact real surface.

From the {\bf birational} point of view, any $\mathbb P^1$-bundle is equivalent
to the trivial bundle: there are infinitely many holomorphic sections
$\sigma:C\to P$; after choosing $3$ distinct ones $\sigma_0$, $\sigma_1$
and $\sigma_\infty$, one defines a birational transformation 
$\phi:P\dashrightarrow C\times\mathbb P^1$ commuting with $\pi$
by sending those sections respectively to $\{y=0\}$, $\{y=1\}$ and $\{y=\infty\}
$.
When the $3$ sections are {\bf disjoint}, the transformation $\phi$
is actually biregular and $P$ is the trivial bundle $C\times\mathbb P^1$.

The {\bf analytic classification} is a much more subtle problem.
If $P$ admits $2$ disjoint sections, say
$\sigma_0,\sigma_\infty:C\to P$, we then say that the bundle is {\bf decomposable}:
one can choose trivialization charts sending those two sections respectively
onto $\{y=0\}$ and $\{y=\infty\}$, so that $P$ may be viewed as
the compactification $\overline L$ of a line bundle $L$. Recall that line bundles
are analytically classified by the Picard group $\text{Pic}(C)$. 
Any two elements $L,L'\in\text{Pic}(C)$
have the same compactification $P$ if, and only if, $L'=L$ or $L^{\otimes(-1)}$:
we just exchange the role of $\sigma_0$ and $\sigma_\infty$ 
(see proof of Proposition \ref{P:RiccatiTorus}).

For instance, on $C=\mathbb P^1$, $\text{Pic}(\mathbb P^1)\simeq\mathbb Z$
and the compactification of $\mathcal O(e)$ (or $\mathcal O(-e)$), 
$e\in\mathbb N$, gives rise to the Hirzebruck surface $\mathbb F_e$.
It follows from Birkhoff's Theorem \cite{Birkhoff} that all 
$\mathbb P^1$-bundle is decomposable on $\mathbb P^1$ and is thus
one of the $\mathbb F_e$ above.

An important analytic invariant of a $\mathbb P^1$-bundle over a curve $C$
is the minimal self-intersection number of a section
$$e(P):=-\text{min}\{\sigma.\sigma\ ;\ \sigma:C\to P\}\in\mathbb Z.$$
For a decomposable bundle $P=\overline L$, $L\in\text{Pic}(C)$, we have 
$e(\overline{L})=\vert\deg(L)\vert\ge 0$. For an undecomposable bundle, Nagata proved 
in \cite{Nagata70} that $-g\le e\le 2g-2$ and all those possibilities occur.

From the {\bf homological} point of view, $H^2(P,\mathbb Z)$ is generated
by the homology class of $\sigma_0$ and $f$ where $\sigma_0$ is any holomorphic 
section
and $f$ any  fibre. Let us choose $\sigma_0$ with minimal self-intersection:
$$\sigma_0\cdot\sigma_0=-e,\ \ \ f\cdot f=0\ \ \ \text{and}\ \ \ \sigma_0\cdot f
=1.$$
The homology class of any other holomorphic section is $\sigma=\sigma_0+n\cdot f
$
with $n\in\mathbb N$: it has self-intersection 
$$\sigma\cdot\sigma=\sigma_0\cdot\sigma_0+2n\cdot\sigma_0\cdot f+n\cdot f\cdot f
=-e+2n\ge -e.$$
In particular, the intersection number of holomorphic sections are either all even, 
either all odd: $e$ mod $2$ is the {\bf topological invariant} of the bundle.

On the other hand, if $\sigma_{0}$ and $\sigma$ are not homologous then
the intersection number $\sigma_0\cdot\sigma=n-e$ must be non
 negative
and we deduce that $\sigma\cdot\sigma\ge e$: when $e>0$, this implies that {\it 
$\sigma_0$
is the unique holomorphic section having negative self-intersection}; there is a
 gap between
$-e$ and $e$. 

\begin{thm}[Atiyah \cite{Atiyah55}]\label{T:Atiyah} 
Beside compactifications of line bundles, there are exactly $2$ undecomposable 
$\mathbb P^1$-bundles over an elliptic curve, $P_0$ and $P_1$, with 
invariant $e=0$ and $-1$ respectively.
\end{thm}

A $\mathbb P^1$-bundle $P$ is {\bf flat} (in the sense of Steenrod \cite{Steenrod51})
when a trivializing atlas can be choosen with constant transition maps 
$\varphi_{i,j}\in\text{PGL}(2,\mathbb C)$ (not depending on $u$).
This means that this atlas defines by the same time a foliation $\mathcal F$
transversal to the fibres on $P$, namely the horizontal foliation
defined by $\{y=\text{constant}\}$ in trivializing coordinates, see Remark~\ref{rm13}.

\begin{thm}[Weil \cite{Weil38}]\label{T:Weil} 
The flat $\mathbb P^1$-bundles over $C$ are all undecomposable bundles
and all those arising as compactification of elements of $\text{Pic}_0(C)$.
\end{thm}

The pairs $(\pi:P\to C,\mathcal F)$ are classified by $H^1(C,\text{PGL}(2,\mathbb C))$.


\subsection{Riccati foliations on $\mathbb P^1$-bundles}\label{sec2.2}

A Riccati foliation on the bundle $\pi:P\to C$ is a singular foliation 
(see definition in  \cite{Brunella})
$\mathcal F$ on $P$ which is transversal to a generic fibre. In trivialization charts $(u_i,y)$,
it is defined by a Riccati differential equation 
${dy\over d{u_i}}=a({u_i})y^2+b({u_i})y+c({u_i})$, 
$a,b,c$ meromorphic in $u$, whence the name. The poles of the coefficients
correspond to vertical invariant fibres for the foliation. Outside of those poles,
the leaves of the foliation are graphs of solutions for the Riccati equation.
The foliation $\mathcal F$
arising in the fibre bundle picture of section \ref{S:BundlePicture} is a {\bf regular} 
Riccati foliation. Nevertheless, we will need to deal with singular foliations
later.

One can define the monodromy representation of a Riccati foliation as
$$\rho:\pi_1(C-\{\text{projection of invariant fibres}\})\to\text{PGL}(2,\mathbb C).$$
A classical Theorem due to Poincar\'e asserts that, in the regular case, the monodromy representation
characterizes the Riccati foliation as well as the $\mathbb P^1$-bundle supporting it
up to analytic equivalence.

\begin{remark}\rm 
One can view a Riccati foliation $\mathcal F$ on the $\mathbb P^1$-bundle
$P=\mathbb PV$ as the projectivization of a meromorphic linear connection
$\nabla$ on the vector bundle $V$. 
In fact, given a (meromorphic linear) connection $\zeta$ on the determinant bundle $\det V=\bigwedge^{2}V\to C$, there is a unique
connection $\nabla$ on $V$ lifting $\mathcal F$ and such that 
$\text{trace}(\nabla)=\zeta$. Indeed, over a local coordinate
$u_i:U_i\to\mathbb C$, the bundle $V$ is trivial and a connection
$\nabla$ is just a meromorphic system
$$\nabla:\ \ \ 
\frac{d}{du_i}\begin{pmatrix}y_1\\ y_2\end{pmatrix}
=\begin{pmatrix}
\alpha(u_{i}) & \beta(u_{i})\\ \gamma(u_{i}) & \delta(u_{i})
\end{pmatrix}
\begin{pmatrix}y_1\\ y_2\end{pmatrix}$$
and the trace of $\nabla$ is the rank $1$ connection defined by 
$$\zeta:=\text{trace}(\nabla):\ \ \ 
\frac{d\lambda}{du_i}
=\left(\alpha(u_{i})+\delta(u_{i})\right)\lambda.$$
The projection of $\nabla$ on $\mathbb PV$ is therefore the Riccati 
equation defined in affine coordinate $(y:1)=(y_1:y_2)$ by
$$\mathcal F:=\mathbb P\nabla:\ \ \ \frac{dy}{du_{i}}=-\gamma(u_{i})y^{2}+(\alpha(u_{i})-\delta(u_{i}))y+\beta(u_{i}).$$
Clearly, $\nabla$ is uniquely defined by $\mathcal F$ and $\zeta$.
Notice finally that the line bundle $\det V$ admits a linear connection $\zeta$ without poles if and only if it belongs to $Pic_{0}(C)$.
\end{remark}

We start recalling some usefull homological formulae from \cite{Brunella}.
First of all, let us introduce $T_\mathcal F$, the {\bf tangent bundle} of $\mathcal F$,
which is a line bundle on the total space $P$ defined as follows.
In trivializing charts $({u_i},y)$, the Riccati foliation is also defined
by the meromorphic vector field 
$V_i:=\partial_{u_i}+(a({u_i})y^2+b({u_i})y+c({u_i}))\partial_y$.
The leaves of $\mathcal F$ are just complex trajectories of the vector field $V_i$.
After choosing a global (meromorphic) vector field $v$ on $C$, one can write 
$v=f_i\cdot\partial_{u_i}$ for a meromorphic function $f_i$ on the chart $U_i$
so that the new meromorphic vector fields $f_i\cdot V_i$ glue together into
a global meromorphic vector field $V$ on $P$ still defining $\mathcal F$
at a generic point. One can think of $V$ as the lifting of $v$ by the (meromorphic)
projective connection defined $\mathcal F$ on the bundle.
Then, $T_\mathcal F$ is the line bundle defined by the divisor of $V$,
i.e. $T_\mathcal F=\mathcal O((V)_0-(V)_\infty)$. If $d$ denotes the number of 
invariant fibre (counted with the multiplicity of the corresponding pole for $V_i$),
then the homology class of $T_\mathcal F$ is given by 
$T_\mathcal F=(2-2g-d)\cdot f$.

Given a curve $\sigma$ on $P$, each component of which is not invariant by $\mathcal F$,
then the number of {\bf tangencies} $\text{Tang}(\mathcal F,\sigma)$ counted with multiplicities
is given by (see \cite{Brunella}, p. 23)
\begin{equation}\label{E:Tang}
\text{Tang}(\mathcal F,\sigma)=\sigma\cdot\sigma-T_\mathcal F\cdot\sigma.
\end{equation}
For instance, if $\sigma=\sigma_0+n\cdot f$ is a section, we immediately deduce that $\text{Tang}(\mathcal F,\sigma)=2n-e-2+2g+d$.

\begin{proof}[Proof of Poincar\'e Theorem \ref{T:Poincare}]
Consider two projective structures on $C$ (compatible with the complex 
structure of $C$) having the same monodromy representation: 
by the construction given in section \ref{S:BundlePicture}, they correspond
to triples $(\pi:P\to C,\mathcal F,\sigma)$ and $(\pi:P\to C,\mathcal F,\sigma')$
with common $\mathbb P^1$-bundle and Riccati foliation.
Since $\mathcal F$ is regular and the section $\sigma$ defining the first projective structure is transversal to $\mathcal F$,
we have $d=0$, $\text{Tang}(\mathcal F,\sigma)=0$, and we deduce that $e=2n+2g-2$.
On the other hand, $\text{Tang}(\mathcal F,\sigma_0)=2g-2-e$ should be non negative
and we obtain $e=2g-2$ and $n=0$: in the genus $g\ge2$ case, $\sigma=\sigma_0$ is the unique 
section having negative self-intersection in $P$, and by the way $\sigma'=\sigma$.
In genus $0$ case, there is nothing to show; in genus $1$ case, the result follows 
directly from formula (\ref{F:LinearMonodromyTorus}) 
and Theorem \ref{T:GunningAffine}.
\end{proof}

Another important formula is the Camacho-Sad Index Theorem (see \cite{Brunella}).
Given a curve $\sigma$ on $P$ invariant by $\mathcal F$, the self-intersection number
of $\sigma$ equals the sum of Camacho-Sad index of $\mathcal F$ along this curve.
When $\mathcal F$ is regular, all invariant curves are smooth and all Camacho-Sad index
vanish: when $\mathcal F$ is regular, any invariant curve $\sigma$ has zero self-intersection.

For instance, if $\mathcal F$ has affine monodromy, then the fixed point gives rise to
an invariant section $\sigma_\infty:C\to P$. We deduce that $e=0$ and $\sigma_\infty$
realizes this minimal self-intersection number. In particular, we recover the fact
that a projective structure on a genus $g\ge2$ curve cannot have affine monodromy
since the corresponding bundle has invariant $e=2g-2>0$. 

More generally, if the monodromy of $\mathcal F$ has a finite orbit (e.g. a finite group
of the infinite dihedral group), then $\mathcal F$ has an invariant curve 
$\sigma=m\cdot\sigma_0+n\cdot f$ and formulae
$\sigma\cdot\sigma=m(2n-em)=0$ together with $\sigma\cdot f=m\ge0$ and 
$\sigma\cdot\sigma_0=n-em\ge0$ show that $e\le0$. Again, this is not
the monodromy of a projective structure whenever $g\ge2$. 

In the particular case where the monodromy of $\mathcal F$ is linear,
we have $2$ invariant disjoint sections $\sigma_0$ and $\sigma_\infty$
showing that the bundle $P$ is actually a compactification of a
line bundle $L\in\text{Pic}_0(C)$.  

We should emphasize that any two
line bundles $L$ and $L'$ have the same compactification if, and only if,
$L'=L$ or $L^{\otimes(-1)}$. Indeed, we first note that, 
for a bundle $P$ satisfying $e(P)=0$, any two sections $\sigma$ and $\sigma'$
are disjoint if, and only if, they have $0$ self-intersection (and are distinct).
The compactification $\overline L$ of a line bundle $L$ always has 
the two canonical disjoint sections $\sigma_0$ and $\sigma_\infty$. 
Now, a diffeomorphism $\phi:\overline L\to\overline L'$ 
between the compactifications of $2$ non trivial line bundles
has to preserve or permute the two canonical sections;
in the former case, $\phi$ is actually an equivalence of line bundles;
in the latter case, $\phi$ restricts to the fibres as $1/z$ and invert
the monodromy. Of course, $\overline L$ is trivial if, and only if, 
$L$ is trivial as a line bundle. It follows that when $C$ has genus $1$ the corresponding set 
of equivalence classes of $\mathbb P^1$-bundles may be thought as
$C/\{\pm1\}\simeq\mathbb P^1$.

A bundle $P$ obtained by suspension of a representation $\rho:\pi_1(C)\to\text{PGL}(2,\mathbb C)$
is topologically trivial ($e$ even) if, and only if, $\rho$ can be lifted as a
representation  $\tilde\rho:\pi_1(C)\to\text{SL}(2,\mathbb C)$.

There is an algebraic and somewhat technical notion of (semi-) stability of vector bundles of arbitrary rank on Riemann surfaces due to Mumford, see for instance \cite{NarasimhanSeshadri65}.
We can define the (semi-) stability of a $\mathbb P^{1}$-bundle $\mathbb PV$ by the same requirement to the rank $2$ vector bundle $V$. It turns out that
a $\mathbb P^1$-bundle is {\bf stable} (resp. {\bf semi-stable}) 
when $e<0$ (resp. $e\le 0$). It is known that if such a bundle occur along an algebraic (resp. analytic) family, it occurs for a Zariski
open subset of the family. 
There is a theorem of Narashimhan and Seshadri characterizing stable bundles on a compact Riemann surface $C$ by means of a precise, but some technical, construction in terms of unitary representations of the fundamental group of $C$. We present here a more comprehensible consequence, see Corollaries 1 and 2 of~\cite{NarasimhanSeshadri65}:

\begin{thm}[Narasimhan-Seshadri \cite{NarasimhanSeshadri65}]
Let $C$ be a compact Riemann surface of genus $g\ge 2$. Then a holomorphic vector bundle of degree zero is stable if and only if it arises from an irreducible unitary representation of the fundamental groups $\pi_{1}(C)$ of $C$. A holomorphic vector bundle on $C$ arises from a unitary representation of the fundamental group if and only if each of its undecomposable components is of degree zero and stable.
\end{thm}

Applying this general result to our situation we obtain that
the map $\rho\mapsto(\pi:P\to C,\mathcal F)\mapsto(\pi:P\to C)$ 
which to a representation $\rho\in\mathrm{Hom}(\pi_1(C),\mathrm{PGL}(2,\mathbb C))$
associate the $\mathbb P^1$-bundles obtained by suspension (forgetting the flat structure)
induces a bijection from the set of irreducible representations $\rho:\pi_1(C)\to\mathrm{PSU}(2,\mathbb C)$ up to $\mathrm{PSU}(2,\mathbb C)$ conjugacy onto the set of isomorphism
class of $\mathbb P^1$-bundles with invariant $e<0$ and even (not fixed).

The complete analytic classification of $\mathbb P^1$-bundles (including unstable ones)
over curves of genus $2$ has been achieved by the works of Atiyah \cite{Atiyah55} and Maruyama
\cite{Maruyama70}. The analytic classification of rank $2$ stable vector bundles over curves of
arbitrary genus from the algebraic point of view (in contrast with Narasimhan-Seshadri's approach)
has been done by Tyurin in \cite{Tyurin64} (see also \cite{Tyurin74} for a survey in arbirtrary rank).

\subsection{Birational geometry of $\mathbb P^1$-bundles}\label{birational}

Given a point $p$ on (the total space of) a $\mathbb P^1$-bundle $\pi:P\to C$,
we will denote by $\text{elm}_pP$ the new $\mathbb P^1$-bundle obtained after 
elementary transformation centered at $p$: after blowing-up the point $p$,
$\text{elm}_pP$ is obtained by contracting the strict transform of the fiber
passing through $p$. 
The strict transform of a section $\sigma$ passing through $p$ (resp. not passing through $p$)
is a section of the new bundle having self-intersection $\sigma\cdot\sigma-1$
(resp. $\sigma\cdot\sigma+1$).
All birational transformations between $\mathbb P^1$-bundles
over curves are obtained by composing finitely many elementary transformations.
On the other hand, any $\mathbb P^1$-bundle over a curve is birational to the
trivial bundle.

\begin{example}\rm 
For instance, let $D$ be a divisor on $C$ and let $p_0$ be the point 
on the zero section of the (total space of the) line bunde $\mathcal O(D)$
over $x\in C$. Then 
$$\text{elm}_{p_0}\overline{\mathcal O(D)}=\overline{\mathcal O(D-[x])}.$$
Similarly, if $p_\infty$ lie on the infinity section of $\overline{\mathcal O(D)}-\mathcal O(D)$ over $x$,
then 
$$\text{elm}_{p_\infty}\overline{\mathcal O(D)}=\overline{\mathcal O(D+[x])}.$$
Now, recall that, as a consequence of Abel Theorem, the map 
$$C^{g}\to\mathrm{Pic}_{0}(C)\ ;\ (x_{1},\ldots,x_{g})\mapsto \mathcal{O}(g[x_{0}]-[x_{1}]-\cdots-[x_{g}])$$
is surjective for any $x_{0}\in C$: it follows that (compactification of) line bundles of degree $0$
can be obtained after applying at most $2g$ elementary transformations to the trivial bundle.
\end{example}

In \cite{MaruyamaNagata69}, Maruyama and Nagata proved that 
an undecomposable bundle can be obtained from the trivial one 
after at most $2g+1$ elementary transformations. On the other hand,
we note that the minimal number of elementary transformations needed
to trivialize all decomposable bundle is unbounded: for a line bundle 
of large degree $d>>0$, one need at least $d$ elementary transformations.

\subsection{Riccati equation, schwartzian derivative and the $2^{\mathrm{nd}}$ order linear differential equation}

First, we would like to make explicit the correspondance between
the point of view of quadratic differentials, and that one of bundle triples.

Consider the triple 
$(\pi:P\to C,\mathcal F,\sigma)$ associated to a projective structure on the curve $C$.
One can reduce $P$ to the trivial bundle and $\sigma$ to the infinity section $\{y=\infty\}$
either locally, by a fibre bundle isomorphism, or globally on $C$, by birational transformation.
Here below, we adopt the later point of view; everything can be carried out {\it mutatis mutandis}
in the local regular setting.
After a birational trivialization like above, $\mathcal F$ becomes possibly singular, but is now defined by a global Riccati equation
\begin{equation}\label{E:GlobalRiccati}
dy+\alpha\cdot y^2+\beta\cdot y+\gamma=0
\end{equation}
where $\alpha,\beta,\gamma$ are meromorphic $1$-forms on $C$. 
This trivialization is unique up to birational transformation of the form
$y=a\tilde y+b$ where $a$ and $b$ are meromorphic function on $C$, $a\not\equiv0$.
Let us see how such change of coordinate acts on the equation.
A change of coordinate of the form $y=a\tilde y$ transforms
the Riccati equation into
\begin{equation}\label{E:ChangeRiccati1}
d\tilde y+a\alpha\tilde y^2+(\beta+{da\over a})\tilde y+{\gamma\over a}=0
\end{equation}
although a change of coordinate $y=\tilde y +b$ yields
\begin{equation}\label{E:ChangeRiccati2}
d\tilde y+\alpha\tilde y^2+(\beta+2b\alpha)\tilde y+(db+b^2\alpha+b\beta+\gamma)=0;
\end{equation}
after a combination of those two transformations, we can choose 
$\alpha$ and $\beta$
arbitrary (with $\alpha\not\equiv0$) and then $\gamma$ is uniquely determined 
by the projective structure.
Let us show how to compute it from the developping map $f$ of the projective structure.

Let us go back to the universal cover where the Riccati foliation is given
by $dy_0=0$ and $\sigma$ is the graph of $f$ (see section \ref{S:BundlePicture}). 
By a preliminary change of coordinate $y_0=y_1+f(u)$, 
we have now $\sigma=\{y_1=0\}$ and the equation becomes $\mathcal F:dy_1+df=0$. 
A second change of coordinate $y_1=f'(u)\cdot y_2$ yields 
$\mathcal F:dy_2+(1+\frac{f''}{f'}y_2)du=0$,
and $\sigma$ is still the zero section $y_2=0$.
In the case of an affine structure on a torus, the later Riccati equation is well-defined:
the corresponding triple $(\pi:P\to C,\mathcal F,\sigma)$ is then given by:
$$P=C\times\mathbb P^1\ni(u,y),\ \ \ \mathcal F:dy+(1+cy)du=0\ \ \ \text{and}\ \ \ 
\sigma(u)\equiv0$$
for some $c\in\mathbb C$.

In the general projective case, it is more convenient to send the section $\sigma$
to the infinity: in the coordinate $\tilde y_2=-1/y_2$, $\mathcal F$ is defined by
$d\tilde y_2+(\tilde y_2^2-\frac{f''}{f'}\tilde y_2)du=0$
and $\sigma$, by $\tilde y_2=\infty$. We finally apply the change of coordinate
$\tilde y_2=y+\frac{1}{2}\frac{f''}{f'}$ and obtain 
\begin{equation}\label{Eq:RiccatiSchwarzian}
\mathcal F: dy+(y^2+\frac{1}{2}\mathcal S_u(f))du=0
\end{equation}
where $\mathcal S_u(f)$ is the schwartzian derivative of $f$
with respect to the variable $u$.
Unfortunately, $du$ is not a global $1$-form. Moreover, 
$u$ is a transcendental variable that we do not want to deal with
when we are considering a triple $(\pi:P\to C,\mathcal F,\sigma)$.
In general, by birational trivialization of the bundle,
one can assume $\sigma$ at infinity and, after choosing a global 
holomorphic $1$-form $\alpha$ on $C$, reduce the Riccati foliation to the special form
\begin{equation}\label{Eq:RiccatiNormalForm}
\mathcal F: d\tilde y+\alpha\tilde y^2+\gamma=0
\end{equation}
with $\gamma$ meromorphic on $C$.
Here, $\alpha$ plays the role of $du$, that is
 $u$ is replaced by a variable $v$ such that $\alpha=dv$;
this takes sense at least at a generic point
of $C$ where everything is regular.
Setting $u=\psi(v)$, the change of coordinate 
$y=\frac{1}{\psi'}\left(\tilde y+\frac{1}{2}\frac{\psi''}{\psi'}\right)$
transforms equation (\ref{Eq:RiccatiSchwarzian}) into 
(\ref{Eq:RiccatiNormalForm}); 
after computation we find
$\gamma=\frac{1}{2}\left( \mathcal S_u(f)\circ\psi\cdot(\psi')^2+\mathcal S_v(\psi) \right)dv$.
Using (\ref{F:SchwartzianComposition}), one finally obtains $\gamma=\mathcal S_v(f\circ\psi)dv$
where $\mathcal S_v$ is the schwarzian derivative with respect to $v$ and deduce

\begin{prop}
Let $(\pi:P\to C,\mathcal F,\sigma)$ be a triple defining a projective structure
on $C$. Let $(v,y)\in U\times\mathbb P^1$ be bundle coordinates over $\pi^{-1}(U)$,  
$U\subset C$, such that
$$\sigma:y=\infty\ \ \ \text{and}\ \ \ \mathcal F:dy+(y^2+\frac{\phi(v)}{2})dv=0.$$
Then the projective coordinates $f$ on $U$ are the solutions of $\mathcal S_vf=\phi$.
\end{prop}

\begin{remark}\rm
Following \cite{MaruyamaNagata69}, the maximally unstable
undecomposable $\mathbb P^1$-bundle $P$ 
corresponding to projective triples $(\mathcal P,\mathcal F,\sigma)$ can be trivialized
after $2g$ elementary transformations (here $e=2g-2$ is even).
The birational tranformation constructed above to put $\mathcal F$
into the normal form  (\ref{Eq:RiccatiNormalForm})
however needs much more elementary transformations.

Indeed, at a point where $\alpha=dv\sim u^\nu du$ has a zero of order $\nu$,
i.e. $v\sim u^{\nu+1}$, the expression 
$$\alpha\otimes\gamma=\frac{1}{2}\mathcal S_v(f)dv^{\otimes 2}\sim\left(\frac{dv}{v}\right)^{\otimes 2}\sim\left(\frac{du}{u}\right)^{\otimes 2}$$
has a pole of order $2$ and thus $\gamma\sim\frac{du}{u^{\nu+2}}$ has a pole of order
$\nu+2$. In fact, $\psi'\sim\frac{1}{u^\nu}$ and the birational change of coordinate
takes the form 
$$y\sim u^\nu\left(\tilde y-\frac{\nu}{2}\frac{1}{u^{\nu+1}}\right)=\frac{1}{u^{\nu+1}}\left(u^{2\nu+1}\tilde y-\frac{\nu}{2}\right);$$
$3\nu+2$ elementary transformations are needed at this point.

Now, we look for a sharp birational trivialization
of $P$, that is to say with exactly $2g$ elementary 
transformations. For any choice of global meromorphic $1$-forms 
$\alpha$ and $\beta$, there is a unique birational transformation
of the form $y=a(\tilde y+b)$ putting the initial Riccati equation
(\ref{Eq:RiccatiSchwarzian}) into the form (\ref{E:GlobalRiccati})
and we have 
$$\left\{\begin{matrix}
\alpha&=&adu\hfill\\
\beta&=&{da\over a}+2ab\hfill\\
\gamma&=&db+ab^2du+{da\over a}b+\frac{\mathcal S_u(f)}{2a}du
\end{matrix}\right.$$
Each zero or pole of $\alpha$ (or $a$)
gives rise to an elementary transformation:
if we choose $\alpha$ holomorphic, we already get
$2g-2$ elementary transformations with the first change
of coordinate. We would like now $ab=\frac{1}{2}(\beta-\frac{da}{a})$ be holomorphic (as much as possible). The sum of residues
of $\frac{da}{a}$ is $2g-2$: we can construct a meromorphic
$1$-form $\beta$ having the same principal part as $\frac{da}{a}$, 
plus one extra simple pole (at, say, $p$) with residue $2-2g$. 
The final change of coordinate $y=a\tilde y+ab$ is therefore a combination of $2g$ elementary transformations: the change of 
coordinate  $y=a\tilde y$ goes from the trivial bundle to 
$\overline{K}$ with $2g-2$ elementary transformations; 
the ultimate transformation $y=\tilde y+ab$ has one simple pole
corresponding to a succession of $2$ generic elementary transformations of the same fibre (compare \cite{Maruyama70}).
\end{remark}

Setting $y=z'/z$, $z'=\frac{dz}{dv}$, the differential equation $dy+(y^2+\frac{\phi(v)}{2})dv=0$
is transformed into 
\begin{equation}\label{E:ScalarEquation}
z''+\frac{\phi(v)}{2}z=0.
\end{equation}
Then the following goes back to Schwarz:

\begin{prop}Any solution $f$ to the differential equation $\mathcal S_v(f)=\phi(v)$
takes the form $f=z_1/z_2$ where $z_1$ and $z_2$ are independant solutions of 
(\ref{E:ScalarEquation}).
\end{prop}

\begin{proof}A straightforward computation shows that $\mathcal S(z_1/z_2)=-2\frac{z_2''}{z_2}$ provided that $z_1$ and $z_2$ are solutions of (\ref{E:ScalarEquation}).
Any other solution $\mathcal S_v(f)=\phi(v)$ takes the form 
$f=\frac{az_1+bz_2}{cz_1+dz_2}$, a quotient of two other solutions.
In fact, one can take $f=z_1/z_2$ with $z_1=\frac{f}{\sqrt{f'}}$ and $z_2=\frac{1}{\sqrt{f'}}$.
\end{proof}

\begin{remark}\rm
One can easily generalize the notion of projective structure to the branching case
by considering triples $(\pi:P\to C,\mathcal F,\sigma)$ with $\sigma$ generically
transversal to $\mathcal F$: branching points of the structure are those points
$x\in C$ over which $\sigma$ has a contact with $\mathcal F$. The local projective
chart then takes the form $f\sim u^{\nu+1}$ where $\nu\in\mathbb N$ is the order 
of contact. More generally, one can consider a singular Riccati foliation $\mathcal F$
generically transversal to $\sigma$, or equivalently linear equation
$dy+(y^2+\frac{\phi(v)}{2})dv=0$ with $\phi$ meromorphic. By the way, projective 
structures on the $3$-punctured sphere ($3$ simple poles) correspond to the Gauss 
Hypergeometric equation, on the $4$-punctured sphere, to the Heun equation and on the punctured torus, to the
Lam\'e equation.
\end{remark}

\begin{remark}\rm
Let $D=\sum\limits_{i=1}^{k}\nu_{i}p_{i}$ be an effective divisor on $\Sigma_{g}$. Consider the set $\mathcal P_{g}(D)$ consisting in all the projective structures on $\Sigma_{g}$ branched over the points $p_{i}$ with ramification order $\nu_{i}\ge 0$ (see \cite{Mandelbaum72,Mandelbaum73}). Notice that the case $D=0$ corresponds to genuine projective structures on $\Sigma_{g}$.
As before we can describe the elements of $\mathcal P_{g}(D)$ as triples $(P,\F,\sigma)$, where $P\to \Sigma_{g}$ is a $\mathbb P^{1}$-bundle with structural group $\PSL$, $\F$ is a transversely projective foliation transverse to the fibres, $\sigma$ is a section such that
$\sigma(p_{i})$ is a tangency point with $\F$ of order $\nu_{i}$ for each $i=1,\ldots,k$ and outside these points $\sigma$ is transverse to $\F$. Projecting to $\Sigma_{g}$ the branched projective structure induced by $\F$ on $\sigma$ we obtain an orbifold complex structure $C$ over $(\Sigma_{g},D)$.
We can make a finite number of elementary transformations 
centered at the tangency points $\sigma(p_{i})$ in order to obtain a birationally equivalent triple $(P',\F',\sigma')$, where $\F'$ is a singular Riccati foliation and $\sigma':C\to P'$ is a holomorphic section everywhere transverse to $\F'$. Applying the same transversality arguments of the proof of Hejal's theorem to $(P',\F',\sigma')$ one shows that the monodromy mapping $\mathcal M:\mathcal P_{g}(D)\to \mathcal R_{g}$ is also a local diffeomorphism.
\end{remark}

\section{The genus $1$ case}

\subsection{Monodromy and bundles}

\begin{prop}\label{P:RiccatiTorus}
Let $C=\mathbb C/(\mathbb Z+\tau\mathbb Z)$ be an elliptic curve, 
$\rho:\pi_1(C)\to\text{PGL}(2,\mathbb C)$ be any representation and
$(\pi:P\to C,\mathcal F)$ be the associated suspension.
Then we are, up to conjugacy, in one of the following cases:
\begin{itemize}
\item $\rho:\pi_1(C)\to\mathbb C^*$ is linear and $P\in\text{Pic}_0(C)$ 
is the compactification of a line bundle; 
$P$ is trivial if, and only if, $\rho(1,\tau)=(e^c,e^{\tau c})$.
\item $\rho:\pi_1(C)\to\mathbb C$ is euclidean and either $P=P_0$ is 
the semi-stable undecomposable bundle, or $P$ is the trivial bundle; 
we are in the latter case if, and only if, $\rho(1,\tau)=(c,\tau c)$.
\item $\rho(1,\tau)=(-z,{1\over z})$ and $P=P_{-1}$ 
is the stable undecomposable bundle.
\end{itemize}
\end{prop}

\begin{proof} 
It is easy to verify that all representation 
$\rho:\mathbb Z^2\to\text{PGL}(2,\mathbb C)$ appear in the statement. 
We have already noticed
that a linear representation gives rise to the compactification of a line bundle (this is almost the definition). In fact, for linear representations, 
we have the exact sequence of sheaves
$$0\to\mathbb C^* \to \mathcal O^*\to \Omega\to 0$$
where $\mathcal O^*$ is the sheaf of invertible holomorphic functions
and the morphism $\mathcal O^*\to \Omega$ is given by $f\mapsto{df\over f}$. From the corresponding exact sequence of cohomology groups, we deduce
the following one
$$0\to H^0(C,\Omega)\to \text{Hom}(\pi_1(C),\mathbb C^*)\to 
\text{Pic}_0(C)\to 0.$$
The first non trivial morphism associates to an holomorphic $1$-form $\omega$
the homomorphism $\gamma\to \exp(\int_\gamma\omega)$ while the second 
one is the suspension.
In our particular case where $C$ is an elliptic curve, we finally deduce
$$0\to \mathbb C du\to \text{Hom}(\pi_1(C),\mathbb C^*)\to C\to 0$$
and the first alternative of the statement follows. 

The suspension of an euclidean representation gives rise to a bundle
with a section $\sigma_\infty$ having $0$ self-intersection.
If there is another section $\sigma_0$ disjoint from $\sigma_\infty$,
then it should be either transversal, or invariant by $\mathcal F$
from (\ref{E:Tang}): in the first case, $\sigma_0$ provides a 
projective structure on $C$ and the monodromy satisfies 
$\rho(1,\tau)=(c,\tau c)$ by Guning Theorem \ref{T:GunningAffine};
in the second case, the monodromy has two fixed points
and is trivial, so is the bundle. In the remaining case
where there is no disjoint section from $\sigma_\infty$,
the bundle is undecomposable with invariant $e=0$ and
we conclude with Atiyah Theorem \ref{T:Atiyah} that $P=P_0$.

Finally, if $\rho(1,\tau)=(-z,{1\over z})$ is the irreducible
representation, we note that $\rho$ cannot be lifted to 
$\text{SL}(2,\mathbb C)$ and thus $e$ is odd. On the other hand,
Weil Theorem \ref{T:Weil} tells us that $P$ must be undecomposable 
(being flat with $e\not=0$). From Atiyah Theorem \ref{T:Atiyah}, 
the only possibility is $P=P_{-1}$.
\end{proof}

\subsection{Algebraic families of bundles and Riccati foliations}

It follows from \cite{Maruyama70} that all degree $0$ line bundles
as well as $P_0$ can be obtained after $2$ elementary transformations
of the trivial bundle. In order to obtain $P_{-1}$, a third one is needed.
We use this approach to provide an algebraic family of flat bundles
and Riccati foliations.

Let $p\in \overline{\mathcal O}$, $q\in \text{elm}_p(C\times\mathbb P^1)$
and consider $P=\text{elm}_q\text{elm}_p\overline{\mathcal O}$.
Fix trivializing coordinates $(u,z)\in C\times\mathbb P^1$ 
and, for simplicity, set $p=(0,\infty)$. This is irrelevant since
all flat $\mathbb P^1$-bundles over $C$ admit a one parameter group of automorphism
lifting the action of $\partial_u$ (see description
of Proposition \ref{P:RiccatiTorus}).
After elementary transformation at the point $p$, we obtain the bundle 
$\overline{\mathcal O(-[0])}$
having one section $\sigma_\infty$ with $-1$ self-intersection and a special
point $\tilde p$, on the fiber over $u=0$ but not on $\sigma_\infty$, 
through which all sections having
$+1$ self-intersection intersect. Indeed, $+1$-sections come from 
horizontal sections of the trivial bundle. Here, we use the fact that there 
is no holomorphic section of homology type $\sigma_0+f$ on $\overline{\mathcal O
}$,
otherwise it would be the graph of a regular covering $C\to\mathbb P^1$.

\begin{figure}[htbp]
\begin{center}

\input{Projective5.pstex_t}
 
\caption{The bundle $\overline{\mathcal O(-[0])}$}
\label{figure:5}
\end{center}
\end{figure}
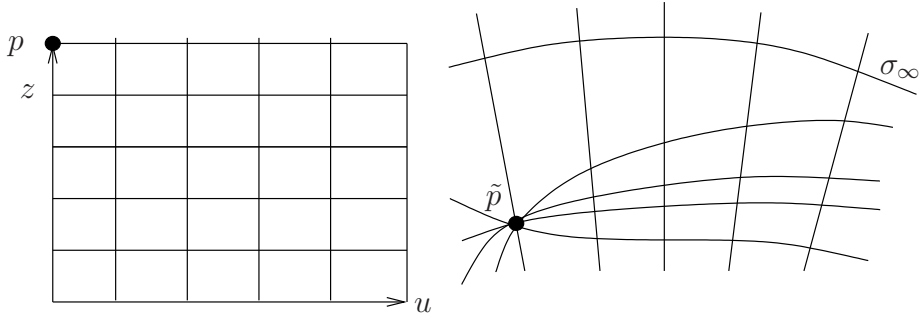

Case 0: $q=\tilde p$. The elementary transformation centered at $\tilde p$ goes 
back 
to the trivial bundle: $P=\overline{\mathcal O}$.

Case 1: $q=(u_0,z_0)$ with $u_0\not=0$ and $z_0\not=\infty$. 
After vertical automorphism, one may assume $z=0$.
The sections $\{z=0\}$ and $\{z=\infty\}$ respectively give rise to
disjoint sections $\sigma_0$ and $\sigma_\infty$ on $P$ 
having $0$ self-intersection. We are in $\text{Pic}_0(C)$ case: $P=\overline L$.

The generic horizontal section $\{z=c\}$ gives rise to a section
$\sigma$ on $P$ intersecting $\sigma_0$ at $u=0$ and $\sigma_\infty$ at $u=u_0$;
in other words, $\sigma$ is a meromorphic section of $L$ with divisor
$\text{Div}(\sigma)=[0]-[u_0]$ on $C$: $L=\mathcal O([u_0]-[0])$
corresponds to $u_0\in\text{Pic}_0(C)\simeq C$.

Case 2: $q$ is on the fibre over $u=0$ but is neither $\tilde p$, nor on
$\sigma_\infty$. Then, $P=P_0$ is the indecomposable bundle.
Indeed, assume that there exists a section $\sigma$ on $P$ disjoint from
$\sigma_\infty$. It then comes from a section of $\overline{\mathcal{O}(-[0])}$ 
disjoint from $\sigma_\infty$
and passing through $q$, itself coming from a section of $\overline{\mathcal O}$
intersecting $\sigma_{\infty}$ only at $u=0$, without multiplicity.
We have already seen that this cannot happen. 

Case 3: $q$ is on $\sigma_{\infty}$, over $u_0$. 
We obtain the bundle $\overline{\mathcal O(-[0]-[u_0])}$.

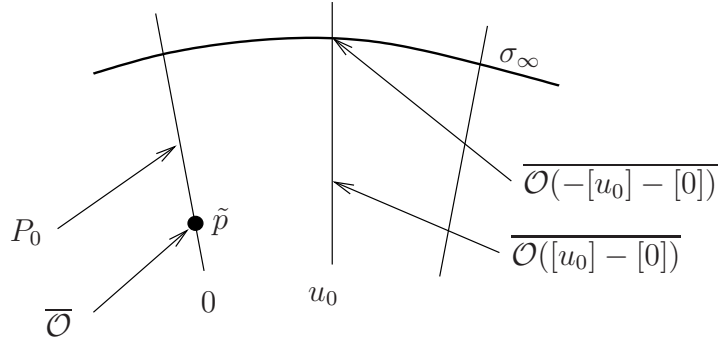
\begin{figure}[htbp]
\begin{center}

\input{Projective6.pstex_t}
 
\caption{An algebraic family of topologically trivial bundles}
\label{figure:6}
\end{center}
\end{figure}

Here, we have parametrized all topologically trivial flat bundles by the line bundle 
$\mathcal O(-[0])$, see \cite{Suwa69}. 
We now want to parametrize all regular Riccati foliations on topologically trivial 
bundles. The natural way to do this is to provide an explicit family of Riccati
equations on the trivial bundle $\overline{\mathcal O}$ having appearant singular 
fibres whose desingularization span all regular Riccati foliations.
For instance, consider a linear Riccati foliation defined on the bundle 
$\overline{\mathcal O([u_0]-[0])}$, $u_0\not=0$. Apart from the invariant sections $\sigma_0$
and $\sigma_\infty$, the leaves are multisections without
zero or pole; after trivialisation of the bundle, those multisections $z(u)$ have 
now a simple pole over $u=0$ and a simple zero over $u=u_0$ (and still have linear
monodromy): their logarithmic derivative ${dz(u)\over z(u)}$ is a meromorphic $1$-form
on $C$ having exactly $2$ simple poles, one at $0$ with residue $-1$ 
and one at $u_0$ with residue $+1$.
In other words, the Riccati equation defining the singular foliation after
trivialization of the bundle is 
$${dz\over du}=\left({\wp'(u)+\wp'(u_0)\over 2(\wp(u)-\wp(u_0))}+c\right)\cdot z
.$$
Indeed, the $1$-form $\left({\wp'(u)+\wp'(u_0)\over 2(\wp(u)-\wp(u_0))}+c\right)
du$
has a simple pole at $u=0$ with residue $-1$ since its principal part is given by ${1\over 2}{\wp'(u)\over \wp(u)}du$ and $\wp$ has a double pole at $u=0$;
the other poles may come from the two zeroes of $\wp(u)-\wp(u_0)$, namely $u=\pm
 u_0$,
but $u=-u_0$ is actually regular since the numerator $\wp'(u)+\wp'(u_0)$
also vanishes at this point: by Residue Theorem, $u=u_0$ is a simple pole with residue $+1$. Of course, any other $1$-form having the same principal part
must differ by a holomorphic $1$-form, namely $c\cdot du$, $c\in\mathbb C$.
We have omited from our discussion the case $u_0=-u_0$ is an order $2$ point 
which can be treated like $u=0$.

After 
two elementary transformations of $\overline{\mathcal{O}}$ centered at
the points $(u,z)=(0,\infty)$ and $(u_0,0)$, 
we obtain by this way all (linear) foliations on the bundle $\mathcal O([u_0]-[0
])$
while $c$ runs over $\mathbb C$. This does not provide
yet a universal family for linear connections on $C$ since the limit
of the Riccati foliation while $u_0\to 0$ is the vertical fibration:
for $(u,z)$ in a compact set not intersecting $\{u=0\}$, $\{z=0\}$ and $\{z=\infty\}$,
we have 
$${\wp'(u)+\wp'(u_0)\over 2(\wp(u)-\wp(u_0))}+c\sim
-{1\over 2}{\wp'(u_0)\over \wp(u_0)}\sim {1\over u_0}\ \ \ \text{while}\ \ \ 
u_0\sim0.$$
In other words, the $1$-form 
$u_0 dz-u_0\left({\wp'(u)+\wp'(u_0)\over 2(\wp(u)-\wp(u_0))}+c\right)\cdot z du$
tends uniformly to $-zdu$ on the compact set, so does the foliation.
We would like to complete this $\mathbb C$-bundle over $u_0\not=0$ with the family
${dz\over du}=c_0\cdot z$, $c_0\in\mathbb C$, of linear connections
on the trivial bundle $\mathcal O$ ($u_0=0$). A way to obtain it from our large 
family 
is obviously to set $c=c(u_0)=c_0-{1\over u_0}$
and take the limit while $u_0\to 0$ with $c_0\in\mathbb C$ fixed.
In other words, in the parameter space $(u_0,c)$ we consider only the limit 
at $(0,\infty)$ while $u_0\to 0$ with a special direction. 
The good global parameter space is obtained after separating
the germs of curves $c+{1\over u_0}=\text{constant}$.
This is done after $2$ elementary transformations on $\overline{\mathcal O}$:
first we blow-up $(u_0,c)=(0,\infty)$ by setting $c=t/u_0$, and then we blow-up
$(u_0,t)=(0,-1)$ by setting $t+1=s u_0$, so that $s={1\over u_0}+c$ coincides 
with the expected parameter $c_0$. The resulting parameter space is the 
{\bf affine bundle}
$A_0:=P_0-\sigma_\infty$ where $\sigma_\infty$ is the unique $0$-section 
of the undecomposable bundle $P_0$.

Now, we construct a fine moduli space as follows.
Consider the product $\overline{\mathcal O}\times\overline{\mathcal O}$
with global coordinates $((u_0,c),(u,z))$, and equipp the bundle over $(u_0,c)$
with the Riccati foliation 
$${dz\over du}=\left({\wp'(u)+\wp'(u_0)\over 2(\wp(u)-\wp(u_0))}+c\right)\cdot z
.$$
This can be seen as an algebraic foliation
on the total space. Now, apply the elementary transformations 
with center along the surfaces $\{u=0,z=\infty\}$ and $\{u=u_0,z=0\}$. 
Then, we modify the base $(u_0,c)\in\overline{\mathcal O}$ 
by two elementary transformations so that we obtain $P_0$ as a base and
the foliation extends as a linear connections all along $u_0=0$.

The euclidean connections on $P_0$ are given by:
$${dz\over du}=\wp(u)+\gamma.$$
Indeed, one can check that the reduction of the singularity over $u=0$
yields $P_0$; on the other hand, it is clear that monodromy is given by translations.
This can be obtained also as a limit of our previous family of connections,
or better from 
$${dz\over du}=\left({\wp'(u)+\wp'(u_0)\over 2(\wp(u)-\wp(u_0))}+c\right)\cdot (
z-c)$$
which is equivalent to the previous one by the change of coordinate $z\mapsto z+
c$.
Now, instead of taking limit along curves $c=c(u_0)=c_0-{1\over u_0}$ as $u_0\to
 0$
with $c_0$ constant, we take limit along $c=c(u_0)=\gamma u_0-{1\over u_0}$,
$\gamma\in\mathbb C$ constant, i.e. $c_0=\gamma u_0$. We have on convenient compact sets:
$${\wp'(u)+\wp'(u_0)\over 2(\wp(u)-\wp(u_0))}-{1\over u_0}\sim
u_0\wp(u)\ \ \ \text{while}\ \ \ 
u_0\sim0$$
so that 
$${dz\over du}=\left({\wp'(u)+\wp'(u_0)\over 2(\wp(u)-\wp(u_0))}+c\right)\cdot (
z-c)\sim \left(\wp(u)+\gamma\right)\cdot (u_0 z+1-\gamma u_0^2)
\sim\wp(u)+\gamma.$$

\subsection{The Riemann-Hilbert Mapping}

For a given elliptic curve $C=\mathbb C/(\mathbb Z+\tau\mathbb Z)$,
the Riemann-Hilbert Mapping provides an analytic isomorphism 
$$\mathcal M:A_0\to\mathbb C^*\times\mathbb C^*$$
between two spaces of algebraic nature. 

The space of linear connections on $C$ is an affine $\mathbb C$-bundle
over $\text{Pic}_0(C)\simeq C$ that we have identified with $A_0$:
it is defined by gluing the chart $(u_0,c)\in (C-\{0\})\times\mathbb C$
with the chart $(u_0,c_0)\in (C,0)\times\mathbb C$ by the transition map
$$(u_0,c)\mapsto(u_0,c_0):=(u_0,c+{1\over u_0}).$$
The space of linear representations of $\pi_1(C)$ is $\mathbb C^*\times\mathbb C^*$.
In the main chart $(u_0,c)$, the analytic connection is given by
$${dz\over z}=\left({\wp'(u)+\wp'(u_0)\over 2(\wp(u)-\wp(u_0))}+c\right)\cdot du.$$
Introducing Weierstrass Zeta Function $\zeta(u)=-\int_0^u\wp(\xi)d\xi$,
one can write $\frac{\wp'(u)+\wp'(u_0)}{2(\wp(u)-\wp(u_0))}
=\zeta(u-u_0)-\zeta(u)+\zeta(u_0)$ and integrate the differential
equation above by means of the Weierstrass Sigma Function: 
the general solution\footnote{This computation was communicated 
to the first author by Frits Beukers; a similar computation 
but with a slightly different presentation was done 
in \cite{Hitchin}.}
is therefore
given by $z(u)=a\frac{\sigma(u-u_0)}{\sigma(u)}e^{\zeta(u_0)\cdot u}$,
$a\in\mathbb C^*$, 
and the monodromy is given by the homomorphism
$$\Lambda=\mathbb Z+\tau\mathbb Z\to\mathbb C^* \ ;\ 
\gamma\mapsto\exp(-u_0\zeta(\gamma)+\zeta(u_0)\gamma+c\gamma).$$
Finally we obtain the full monodromy mapping
$$\mathcal M:A_0\to\mathbb C^*\times\mathbb C^*\ ;\ 
\left\{\begin{array}{rcl}(u_0,c)&\mapsto&
(e^{-u_0\zeta(1)+\zeta(u_0)+c},e^{-u_0\zeta(\tau)+\zeta(u_0)\tau+c\tau})\\
(0,c_0)&\mapsto&(e^{c_0},e^{c_0\tau})\end{array}\right.$$

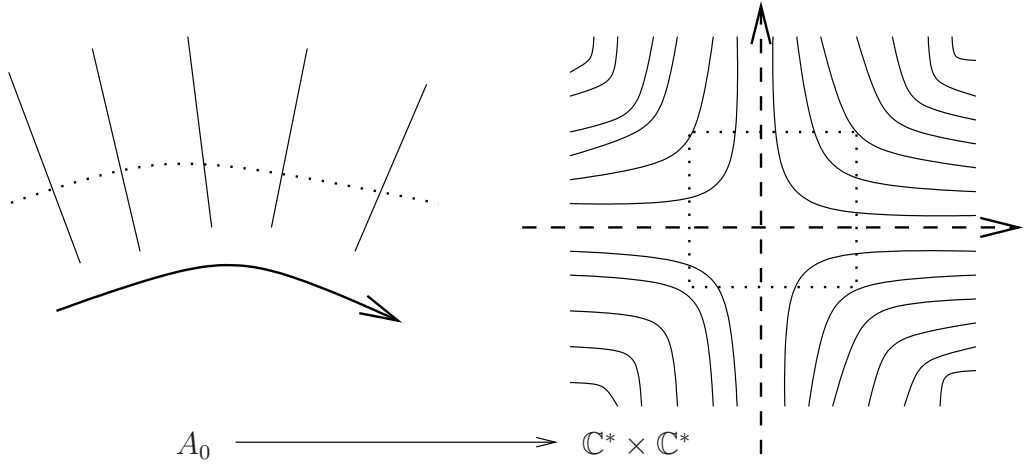
\begin{figure}[htbp]
\begin{center}

\input{Projective7.pstex_t}
 
\caption{The Riemann-Hilbert Mapping}
\label{figure:7}
\end{center}
\end{figure}

The image by the monodromy map of the algebraic fibration defined on $A_0$ 
is the holomorphic foliation defined on $\mathbb C^*\times\mathbb C^*$
by the linear vector field
$$x\partial_{x}+\tau y\partial_{y}.$$

As a particular case of Narasimhan-Seshadri Theorem, the unitary representations
$\mathbb S^1\times\mathbb S^1$ form a smooth real $2$-dimensional torus
transversal to the foliation and cutting each leaf once. It is the space of the 
leaves.
It inherits, from the transversal complex foliation, a complex structure,
namely the structure of $C$. 
The euclidean foliations defined on $A_0$ by ${dc\over du_0}=\wp(u_0)+\gamma$,
$\gamma\in\mathbb C$, are sent to the linear foliations $x\partial x+\lambda y\partial y$,
$\lambda\in\mathbb P^1\setminus\{\tau\}$. 
The space of linear connections is equipped with a group law given by tensor
product; it is just the pull-back of the natural group law 
on $\mathbb C^*\times\mathbb C^*$. We thus get an analytic isomorphism 
between two algebraic groups that are not algebraically equivalent.

One can compute the group law on $A_0$ as follows. Given two connections
$${dz\over z}=\left({\wp'(u)+\wp'(u_1)\over 2(\wp(u)-\wp(u_1))}+c_1\right)\cdot 
du\ \ \
\text{and}\ \ \ 
{dz\over z}=\left({\wp'(u)+\wp'(u_2)\over 2(\wp(u)-\wp(u_2))}+c_2\right)\cdot du
,$$
the tensor product is a connection of the form
$${dz\over z}=\left({\wp'(u)+\wp'(u_3)\over 2(\wp(u)-\wp(u_3))}+c_3\right)\cdot 
du$$
with $u_3=u_1+u_2$ (group law on $\text{Pic}_0(C)$). Then, $c_3$ is determined by the fact
that  
$${\wp'(u)+\wp'(u_1)\over 2(\wp(u)-\wp(u_1))}+{\wp'(u)+\wp'(u_2)\over 2(\wp(u)-\
wp(u_2))}
-{\wp'(u)+\wp'(u_3)\over 2(\wp(u)-\wp(u_3))}+c_1+c_2-c_3={df\over f}$$
for a meromorphic function $f$ on $C$. Looking at the principal part of the left
 hand side,
one see that $f$ must have divisor $\text{Div}(f)=[u_1]+[u_2]-[u_3]-[0]$ so that,
up to a scalar, we have
$$f={\wp'(u)-\wp'(u_1)-{\wp'(u_2)-\wp'(u_1)\over\wp(u_2)-\wp(u-1)}(\wp(u)-\wp(u_
1))
\over\wp(u)-\wp(u_3)};$$
after computations, one finds that
$$c_3=c_1+c_2-{\wp'(u_2)-\wp'(u_1)\over 2(\wp(u_2)-\wp(u_1))}.$$

\end{document}

%% file: Projective1.pstex_t
\begin{picture}(0,0)%
\includegraphics{Projective1.pstex}%
\end{picture}%
\setlength{\unitlength}{3947sp}%
\begingroup\makeatletter\ifx\SetFigFont\undefined%
\gdef\SetFigFont#1#2#3#4#5{%
  \reset@font\fontsize{#1}{#2pt}%
  \fontfamily{#3}\fontseries{#4}\fontshape{#5}%
  \selectfont}%
\fi\endgroup%
\begin{picture}(3965,3131)(1642,-4760)
\put(2551,-3661){\makebox(0,0)[lb]{\smash{\SetFigFont{12}{14.4}{\familydefault}{\mddefault}{\updefault}$f_i$}}}
\put(4351,-3511){\makebox(0,0)[lb]{\smash{\SetFigFont{12}{14.4}{\familydefault}{\mddefault}{\updefault}$f_j$}}}
\put(3226,-2986){\makebox(0,0)[lb]{\smash{\SetFigFont{12}{14.4}{\familydefault}{\mddefault}{\updefault}$U_i$}}}
\put(3976,-2386){\makebox(0,0)[lb]{\smash{\SetFigFont{12}{14.4}{\familydefault}{\mddefault}{\updefault}$U_j$}}}
\put(3526,-4186){\makebox(0,0)[lb]{\smash{\SetFigFont{12}{14.4}{\familydefault}{\mddefault}{\updefault}$\varphi_{i,j}$}}}
\put(5401,-3211){\makebox(0,0)[lb]{\smash{\SetFigFont{12}{14.4}{\familydefault}{\mddefault}{\updefault}$\Sigma_g$}}}
\put(4201,-4711){\makebox(0,0)[lb]{\smash{\SetFigFont{12}{14.4}{\familydefault}{\mddefault}{\updefault}$\mathbb P^1$}}}
\end{picture}

%% file: Projective2.pstex_t
\begin{picture}(0,0)%
\includegraphics{Projective2.pstex}%
\end{picture}%
\setlength{\unitlength}{3947sp}%
\begingroup\makeatletter\ifx\SetFigFont\undefined%
\gdef\SetFigFont#1#2#3#4#5{%
  \reset@font\fontsize{#1}{#2pt}%
  \fontfamily{#3}\fontseries{#4}\fontshape{#5}%
  \selectfont}%
\fi\endgroup%
\begin{picture}(5797,3011)(364,-3671)
\put(1501,-1261){\makebox(0,0)[lb]{\smash{\SetFigFont{12}{14.4}{\familydefault}{\mddefault}{\updefault}$\varphi_1$}}}
\put(1501,-1936){\makebox(0,0)[lb]{\smash{\SetFigFont{12}{14.4}{\familydefault}{\mddefault}{\updefault}$\varphi_2$}}}
\put(1576,-2686){\makebox(0,0)[lb]{\smash{\SetFigFont{12}{14.4}{\familydefault}{\mddefault}{\updefault}$\varphi_3$}}}
\put(676,-1561){\makebox(0,0)[lb]{\smash{\SetFigFont{12}{14.4}{\familydefault}{\mddefault}{\updefault}$\Delta_1^ -$}}}
\put(676,-2236){\makebox(0,0)[lb]{\smash{\SetFigFont{12}{14.4}{\familydefault}{\mddefault}{\updefault}$\Delta_2^-$}}}
\put(676,-2911){\makebox(0,0)[lb]{\smash{\SetFigFont{12}{14.4}{\familydefault}{\mddefault}{\updefault}$\Delta_3^-$}}}
\put(2476,-1186){\makebox(0,0)[lb]{\smash{\SetFigFont{12}{14.4}{\familydefault}{\mddefault}{\updefault}$\Delta_1^+$}}}
\put(2551,-1861){\makebox(0,0)[lb]{\smash{\SetFigFont{12}{14.4}{\familydefault}{\mddefault}{\updefault}$\Delta_2^+$}}}
\put(2401,-2986){\makebox(0,0)[lb]{\smash{\SetFigFont{12}{14.4}{\familydefault}{\mddefault}{\updefault}$\Delta_3^+$}}}
\put(2851,-3511){\makebox(0,0)[lb]{\smash{\SetFigFont{12}{14.4}{\familydefault}{\mddefault}{\updefault}$\mathbb P^1$}}}
\put(5776,-3586){\makebox(0,0)[lb]{\smash{\SetFigFont{12}{14.4}{\familydefault}{\mddefault}{\updefault}$C$}}}
\end{picture}

%% file: Projective3bis.pstex_t
\begin{picture}(0,0)%
\includegraphics{Projective3bis.pstex}%
\end{picture}%
\setlength{\unitlength}{3947sp}%
\begingroup\makeatletter\ifx\SetFigFont\undefined%
\gdef\SetFigFont#1#2#3#4#5{%
  \reset@font\fontsize{#1}{#2pt}%
  \fontfamily{#3}\fontseries{#4}\fontshape{#5}%
  \selectfont}%
\fi\endgroup%
\begin{picture}(5982,2478)(559,-3580)
\put(601,-1261){\makebox(0,0)[lb]{\smash{{\SetFigFont{12}{14.4}{\familydefault}{\mddefault}{\updefault}$\mathbb P^ 1$}}}}
\put(3076,-3436){\makebox(0,0)[lb]{\smash{{\SetFigFont{12}{14.4}{\familydefault}{\mddefault}{\updefault}$U$}}}}
\put(6526,-1486){\makebox(0,0)[lb]{\smash{{\SetFigFont{12}{14.4}{\familydefault}{\mddefault}{\updefault}$P$}}}}
\put(6376,-2236){\makebox(0,0)[lb]{\smash{{\SetFigFont{12}{14.4}{\familydefault}{\mddefault}{\updefault}$\mathcal F$}}}}
\put(6226,-2761){\makebox(0,0)[lb]{\smash{{\SetFigFont{12}{14.4}{\familydefault}{\mddefault}{\updefault}$\sigma$}}}}
\put(6076,-3511){\makebox(0,0)[lb]{\smash{{\SetFigFont{12}{14.4}{\familydefault}{\mddefault}{\updefault}$C$}}}}
\put(2776,-1336){\makebox(0,0)[lb]{\smash{{\SetFigFont{12}{14.4}{\familydefault}{\mddefault}{\updefault}$\text{graph}(f)$}}}}
\end{picture}%

%% file: Projective4.pstex_t
\begin{picture}(0,0)%
\includegraphics{Projective4.pstex}%
\end{picture}%
\setlength{\unitlength}{3947sp}%
\begingroup\makeatletter\ifx\SetFigFont\undefined%
\gdef\SetFigFont#1#2#3#4#5{%
  \reset@font\fontsize{#1}{#2pt}%
  \fontfamily{#3}\fontseries{#4}\fontshape{#5}%
  \selectfont}%
\fi\endgroup%
\begin{picture}(4642,3128)(3064,-4351)
\end{picture}%

%% file: Projective4bis.pstex_t
\begin{picture}(0,0)%
\includegraphics{Projective4bis.pstex}%
\end{picture}%
\setlength{\unitlength}{3947sp}%
\begingroup\makeatletter\ifx\SetFigFont\undefined%
\gdef\SetFigFont#1#2#3#4#5{%
  \reset@font\fontsize{#1}{#2pt}%
  \fontfamily{#3}\fontseries{#4}\fontshape{#5}%
  \selectfont}%
\fi\endgroup%
\begin{picture}(4664,3128)(3064,-4351)
\end{picture}%

%% file: Projective5.pstex_t
\begin{picture}(0,0)%
\includegraphics{Projective5.pstex}%
\end{picture}%
\setlength{\unitlength}{3947sp}%
\begingroup\makeatletter\ifx\SetFigFont\undefined%
\gdef\SetFigFont#1#2#3#4#5{%
  \reset@font\fontsize{#1}{#2pt}%
  \fontfamily{#3}\fontseries{#4}\fontshape{#5}%
  \selectfont}%
\fi\endgroup%
\begin{picture}(5717,2011)(376,-2210)
\put(451,-811){\makebox(0,0)[lb]{\smash{\SetFigFont{12}{14.4}{\familydefault}{\mddefault}{\updefault}$z$}}}
\put(3376,-1486){\makebox(0,0)[lb]{\smash{\SetFigFont{12}{14.4}{\familydefault}{\mddefault}{\updefault}$\tilde p$}}}
\put(376,-511){\makebox(0,0)[lb]{\smash{\SetFigFont{12}{14.4}{\familydefault}{\mddefault}{\updefault}$p$}}}
\put(2926,-2161){\makebox(0,0)[lb]{\smash{\SetFigFont{12}{14.4}{\familydefault}{\mddefault}{\updefault}$u$}}}
\put(5851,-661){\makebox(0,0)[lb]{\smash{\SetFigFont{12}{14.4}{\familydefault}{\mddefault}{\updefault}$\sigma_\infty$}}}
\end{picture}

%% file: Projective6.pstex_t
\begin{picture}(0,0)%
\includegraphics{Projective6.pstex}%
\end{picture}%
\setlength{\unitlength}{3947sp}%
\begingroup\makeatletter\ifx\SetFigFont\undefined%
\gdef\SetFigFont#1#2#3#4#5{%
  \reset@font\fontsize{#1}{#2pt}%
  \fontfamily{#3}\fontseries{#4}\fontshape{#5}%
  \selectfont}%
\fi\endgroup%
\begin{picture}(3472,2170)(2401,-2369)
\put(3676,-1636){\makebox(0,0)[lb]{\smash{\SetFigFont{12}{14.4}{\familydefault}{\mddefault}{\updefault}$\tilde p$}}}
\put(5476,-586){\makebox(0,0)[lb]{\smash{\SetFigFont{12}{14.4}{\familydefault}{\mddefault}{\updefault}$\sigma_\infty$}}}
\put(2401,-1711){\makebox(0,0)[lb]{\smash{\SetFigFont{12}{14.4}{\familydefault}{\mddefault}{\updefault}$P_0$}}}
\put(3601,-2161){\makebox(0,0)[lb]{\smash{\SetFigFont{12}{14.4}{\familydefault}{\mddefault}{\updefault}$0$}}}
\put(4276,-2086){\makebox(0,0)[lb]{\smash{\SetFigFont{12}{14.4}{\familydefault}{\mddefault}{\updefault}$u_0$}}}
\put(2626,-2311){\makebox(0,0)[lb]{\smash{\SetFigFont{12}{14.4}{\familydefault}{\mddefault}{\updefault}$\overline{\mathcal O}$}}}
\put(5551,-1861){\makebox(0,0)[lb]{\smash{\SetFigFont{12}{14.4}{\familydefault}{\mddefault}{\updefault}$\overline{\mathcal O([u_0]-[0])}$}}}
\put(5626,-1411){\makebox(0,0)[lb]{\smash{\SetFigFont{12}{14.4}{\familydefault}{\mddefault}{\updefault}$\overline{\mathcal O(-[u_0]-[0])}$}}}
\end{picture}

%% file: Projective7.pstex_t
\begin{picture}(0,0)%
\includegraphics{Projective7.pstex}%
\end{picture}%
\setlength{\unitlength}{3947sp}%
\begingroup\makeatletter\ifx\SetFigFont\undefined%
\gdef\SetFigFont#1#2#3#4#5{%
  \reset@font\fontsize{#1}{#2pt}%
  \fontfamily{#3}\fontseries{#4}\fontshape{#5}%
  \selectfont}%
\fi\endgroup%
\begin{picture}(6419,2930)(54,-2894)
\put(3676,-2836){\makebox(0,0)[lb]{\smash{\SetFigFont{12}{14.4}{\familydefault}{\mddefault}{\updefault}$\mathbb C^*\times\mathbb C^*$}}}
\put(1126,-2836){\makebox(0,0)[lb]{\smash{\SetFigFont{12}{14.4}{\familydefault}{\mddefault}{\updefault}$A_0$}}}
\end{picture}